\documentclass[12pt,oneside,reqno]{amsart}
\usepackage{graphicx}
\usepackage{mathrsfs}
\usepackage{stmaryrd}
\usepackage{amsfonts}
\usepackage{cite}
\usepackage{enumerate,amsmath,amssymb,amsthm}
\usepackage{booktabs} 
\usepackage{diagbox} 

\newcommand{\dif}{\mathrm{d}}

\newcommand{\be}{\begin{eqnarray}}
	\newcommand{\ee}{\end{eqnarray}}
\newcommand{\ce}{\begin{eqnarray*}}
	\newcommand{\de}{\end{eqnarray*}}
\newtheorem{theorem}{Theorem}[section]
\newtheorem{lemma}[theorem]{Lemma}
\newtheorem{remark}[theorem]{Remark}
\newtheorem{definition}[theorem]{Definition}
\newtheorem{proposition}[theorem]{Proposition}
\newtheorem{Example}[theorem]{Example}
\newtheorem{corollary}[theorem]{Corollary}
\newtheorem{condition}[theorem]{Condition}
\def\e{\varepsilon}
\def\t{\theta}
\def\a{\alpha}

\def\b{\beta}
\def\d{\delta}

\def\g{\gamma}
\def\s{\sigma}

\def\l{\lambda}

\def\[{{\Big[}}
\def\]{{\Big]}}
\def\<{{\langle}}
\def\>{{\rangle}}
\def\({{\Big(}}
\def\){{\Big)}}

\def\no{\nonumber}
\def\bt{\begin{theorem}}
	\def\et{\end{theorem}}
\def\bl{\begin{lemma}}
	\def\el{\end{lemma}}
\def\br{\begin{remark}}
	\def\er{\end{remark}}
\def\bx{\begin{Example}}
	\def\ex{\end{Example}}
\def\bd{\begin{definition}}
	\def\ed{\end{definition}}
\def\bp{\begin{proposition}}
	\def\ep{\end{proposition}}
\def\bc{\begin{corollary}}
	\def\ec{\end{corollary}}
\def\bco{\begin{condition}}
	\def\eco{\end{condition}}
\def\cA{{\mathcal A}}
\def\cB{{\mathcal B}}

\def\cD{{\mathcal D}}

\def\cG{{\mathcal G}}

\def\cI{{\mathcal I}}

\def\cM{{\mathcal M}}

\def\mE{{\mathbb E}}

\def\mM{{\mathbb M}}
\def\mN{{\mathbb N}}

\def\mP{{\mathbb P}}

\def\mR{{\mathbb R}}

\def\mU{{\mathbb U}}
\def\mV{{\mathbb V}}
\def\mW{{\mathbb W}}
\def\mX{{\mathbb X}}

\def\sA{{\mathscr A}}
\def\sB{{\mathscr B}}

\def\sF{{\mathscr F}}
\def\sG{{\mathscr G}}

\def\sL{{\mathscr L}}

\def\sV{{\mathscr V}}

\def\geq{\geqslant}
\def\leq{\leqslant}

\begin{document}
	
\allowdisplaybreaks
\title{Large deviations for invariant measures of multivalued stochastic differential equations with jumps}
	
\author{Huijie Qiao}

	
	
\thanks{This work was supported by NSF of China (No.12071071) and the Jiangsu Provincial Scientific Research Center of Applied Mathematics (No. BK20233002).}
	
	
\date{}
	
\dedicatory{School of Mathematics,
		Southeast University\\
		Nanjing, Jiangsu 211189, China\\
		hjqiaogean@seu.edu.cn}
	
\begin{abstract}
This work focuses on multivalued stochastic differential equations with jumps. First, by employing the weak convergence approach, we establish the Freidlin-Wentzell uniform large deviation principle and the Dembo-Zeitouni uniform large deviation principle for these equations. Subsequently, based on these results, we derive both upper and lower bounds for the large deviations of invariant measures associated with the equations.
\end{abstract}
	
\maketitle \rm
	
\section{Introduction}

Consider the following multivalued stochastic differential equation (SDE for short) with jumps on $\mR^d$:
\be\left\{\begin{array}{ll}
\dif X^\e_t\in -A(X^\e_t)\dif t+b(X^\e_t)\dif t+\sqrt \e\sigma(X^\e_t)\dif W_t+\e\int_{\mU}f(X^\e_{t-},u) \tilde N^{\e^{-1}}(\dif t\dif u),\\
X^\e_0=x_0\in \overline{\cD(A)}, \quad t\geq 0,
\end{array}
\right.
\label{msdej1}
\ee
where $A$ is a maximal monotone operator on $\mR^d$ (See Subsection \ref{mmo}), the coefficients $b: \mR^d\rightarrow \mR^d, \sigma: \mR^d\rightarrow\mR^{d\times l}, f: \mR^d\times\mU\rightarrow\mR^d$ are Borel 
measurable, $\mU$ is a locally compact Polish space and $0<\e<1$ is a small parameter. Here $(W_t)$ is an $l$-dimensional Brownian motion, $N^{\e^{-1}}(\dif t\dif u)$ is a Poisson random measure with the intensity 
$\e^{-1}\dif t\nu(\dif u)$, where $\nu$ is a $\sigma$-finite measure on $\mU$, $\tilde N^{\e^{-1}}(\dif t\dif u):=N^{\e^{-1}}(\dif t\dif u)-\e^{-1}\dif t\nu(\dif u)$ is the compensated martingale measure of $N^{\e^{-1}}(\dif t\dif u)$, 
which are all defined on the complete filtered probability space $(\Omega, \mathscr{F}, \{\mathscr{F}_t\}_{t\in[0,T]}, \mP)$. Moreover, $(W_t)$ and $N^{\e^{-1}}(\dif t\dif u)$ are mutually independent. Equations of the type represented by Eq.(\ref{msdej1}) commonly arise in physics, biology, chemistry, and other scientific disciplines, where they are used to model various constrained dynamical systems (\cite{ff, hjh, lw}). In the case where $A=0$, Eq.(\ref{msdej1}) reduces to a standard SDE with jumps. Consequently, the results established for Eq.(\ref{msdej1}) are applicable to general SDEs with jumps. Furthermore, there exist several existing studies on the well-posedness and ergodicity of Eq.(\ref{msdej1}) (\cite{gw, mr}).

The theory of large deviation principles (LDPs), developed in the 1960s by Freidlin, Wentzell, and others, characterizes the asymptotic decay rates of probabilities of rare events (\cite{q1, rxz, rwz}). Uniform LDPs, in particular, are instrumental in determining the exit times and locations of stochastic processes from given domains (\cite{dz, fw}). A substantial body of research has established uniform LDP results for SDEs (\cite{by, jwzz, Mar, ms, wy, jw}) and stochastic partial differential equations (SPDEs) (\cite{bdv1, cm, cr2, eg, km, zq, ms1, sbd, ss, wang1, wang2}). We now highlight several results closely related to our work. Bao and Yuan \cite{by} established the Freidlin-Wentzell uniform LDP for neutral functional SDEs with jumps. In infinite-dimensional settings, Maroulas \cite{Mar} obtained a Freidlin-Wentzell uniform LDP for Eq.(\ref{msdej1}) involving only the coefficient $f$. Furthermore, Wu \cite{jw} assumed that the coefficients $b$, $\sigma$, and $f$ are bounded and that $f$ depends implicitly on $u$, under which conditions the uniform Laplace principle was derived for Eq.(\ref{msdej1}) with operators $A$ and $f$. However, these assumptions exclude several important classes of equations. Therefore, the first objectives of this paper is to relax certain restrictions imposed in \cite{jw} and establish the Freidlin-Wentzell uniform LDP for Eq.(\ref{msdej1}) under more general conditions. 

In establishing the Freidlin-Wentzell uniform LDP for Eq.(\ref{msdej1}), we employ the weak convergence approach. However, verifying the two sufficient conditions (see Condition \ref{cond}), particularly the first one, is nontrivial. Since Eq.(\ref{msdej1}) involves a maximal monotone operator $A$, the corresponding skeleton equation (\ref{skeleq}) also depends on $A$. As a result, direct estimation of the increment $X^{\e,\xi^\e,\zeta^\e}_{\tau+\d}(x_0)-X^{\e,\xi^\e,\zeta^\e}_{\tau}(x_0)$, where $\tau$ is a stopping time and $(X^{\e,\xi^\e,\zeta^\e}(x_0),K^{\e,\xi^\e,\zeta^\e}(x_0))$ denotes the solution to Eq.(\ref{skeleq}), is not feasible. Consequently, the classical tightness-based argument for $X^{\e,\xi^\e,\zeta^\e}(x_0)$ as employed in \cite{by} does not apply. To overcome this difficulty, we adopt an alternative strategy: first, we establish that $X^{\e,\xi^\e, \zeta^\e}(x_0^\e)$ converges in probability to $X^{\xi,\zeta}(x_0)$ as $x_0^\e \to x_0$ and $(\xi^\e, \zeta^\e)$ converges almost surely to $(\xi, \zeta)$; second, we invoke the Skorohod representation theorem to complete the verification of the first condition. This method has been previously utilized in \cite{jw}; however, therein the author relies on the uniform convergence topology, which is {\it not} suitable for the solution process $X^{\e,\xi^\e, \zeta^\e}(x_0^\e)$ due to the c\`adl\`ag nature of its sample paths in time $t$. In contrast, we employ the Skorohod topology to characterize convergence in $D([0,T],\overline{\cD(A)})$, the space of right-continuous functions with left limits defined on $[0,T]$ and taking values in $\overline{\cD(A)}$. This choice necessitates more intricate technical arguments and renders several proofs significantly more involved.

Next, to investigate large deviations for invariant measures associated with Eq.(\ref{msdej1}), we also establish the Dembo-Zeitouni uniform LDP for this equation. While it is well known that the Freidlin-Wentzell LDP is equivalent to the Dembo-Zeitouni LDP when the rate function is good, such equivalence does not extend to their uniform counterparts (\cite{ms}). Therefore, we prove the continuity of the level sets of the rate function in the Hausdorff metric to obtain the Dembo-Zeitouni uniform LDP for Eq.(\ref{msdej1}).

Under suitable assumptions on $b$, $\s$, and $f$, Eq.(\ref{msdej1}) admits an invariant probability measure $\mu^\e$ (\cite{gw}). A natural question arises: as $\e \to 0$, does $\mu^\e$ converge weakly to $\mu^0$, where $\mu^0$ denotes an invariant probability measure of the following multivalued differential equation?
\be\left\{\begin{array}{ll}
\dif X^0_t\in -A(X^0_t)\dif t+b(X^0_t)\dif t, \\
X^0_0=x_0\in \overline{\cD(A)}, \quad t\geq 0.
\end{array}
\right.
\label{x0eq}
\ee
Furthermore, if such convergence holds, what is the corresponding exponential convergence rate? These questions are closely related to the stability properties of solutions to Eq.(\ref{msdej1}) (\cite{cdjz, hjly}). For instances of Eq.(\ref{msdej1}) without the operators $A$ and $f$, extensive results exist regarding these issues--see \cite{fw, hjly} for finite-dimensional state spaces and \cite{bc, bfz, cp, cr1, dm1, dm2, so, wangiv1, wangiv2, zl} for infinite-dimensional settings.

For Eq.(\ref{msdej1}) involving only the jump term $f$, Ma and Xi \cite{xmfx} studied these problems using a large deviation approach. Subsequently, Chen et al. \cite{cdjz} established convergence of invariant measures for $2D$ stochastic Navier-Stokes equations driven by L\'evy noise via tightness arguments. For Eq.(\ref{msdej1}) involving only the maximal monotone operator $A$, Zhang \cite{hz} derived the LDP for invariant measures of multivalued SDEs. In infinite-dimensional settings, Zhang \cite{tz} proved large deviations for invariant measures of SPDEs with two reflecting walls. To the best of our knowledge, no existing work addresses large deviations for invariant measures of Eq.(\ref{msdej1}) in the presence of both $A$ and $f$. Hence, our second objective is to study the LDP for the invariant measures of Eq.(\ref{msdej1}), specifically by establishing the exponential convergence rate of $\mu^\e$ to $\mu^0$. 

In proving the LDP for the invariant measures of Eq.(\ref{msdej1}), the appearance of the multivalued operator $A$ presents two major challenges. The first challenge stems from the inability to directly obtain moment estimates for the solution of Eq.(\ref{x0eq}), which are crucial for proving the compactness of the level sets of the rate function $V$. To overcome this, we construct an approximating equation using the Yosida approximation of $A$, and then derive the required estimate by analyzing the corresponding solution of the approximate equation. The second challenge lies in establishing the exponential tightness of the solution to Eq.(\ref{msdej1}), a property essential for proving the exponential tightness of the invariant measures. On one hand, the application of the Freidlin-Wentzell uniform LDP (cf. \cite{xmfx,hz}) is not feasible, as it requires handling arbitrarily large initial values. On the other hand, the standard approach of decomposing the solution into four components and estimating each term separately (cf. \cite{tz, so}) cannot be applied here due to the lack of control over the term $|K^\e_t(x_0)|$. To resolve this issue, we establish an exponential moment estimate for the solution of Eq.(\ref{msdej1}) (Lemma \ref{expoesti}) and employ this estimate together with Fatou's lemma to prove the exponential tightness of the invariant measures. Finally, we highlight two key technical aspects of our approach. First, instead of proving the lower semicontinuity of $V$, which is difficult in this setting, we use an equivalent characterization to verify the closedness of its level sets. Second, by exploiting the definition of $V$, we directly establish an inclusion relation among its level sets (Lemma \ref{kbrks}), thus avoiding the need to establish continuity of $V$.

In summary, the novelty of this paper lies in two aspects. Firstly, we prove the Freidlin-Wentzell uniform LDP and the Dembo-Zeitouni uniform LDP for Eq.(\ref{msdej1}). Our Freidlin-Wentzell uniform LDP result encompasses Theorem 3.1 in \cite{jw} and Theorem 2.2 in \cite{by} to some extent. Secondly, large deviations for invariant measures of Eq.(\ref{msdej1}) are established. This result can cover \cite[Theorem 3.2]{xmfx} and \cite[Theorem 3.4]{hz}.

The rest of this paper is organized as follows. In Section \ref{pre}, we introduce some notations, maximal monotone operators and uniform LDPs. The main results are  formulated in Section \ref{main}. And the proofs of main results are placed in Section \ref{uldpthproo} and \ref{ldpinvameasproo}, respectively.
	
The following convention will be used throughout the paper: $C$ with or without indices will denote different positive constants whose values may change from one place to another.

\section{Preliminaries}\label{pre}

In this section, we introduce some notations, maximal monotone operators and uniform LDPs.

\subsection{Notation}

In this subsection, we introduce some notations used in the sequel.

Let $|\cdot|, \|\cdot\|$ be the norms of a vector and a matrix, respectively. Let $\langle\cdot,\cdot\rangle$ be the inner product of vectors on $\mR^d$. $U^{*}$ denotes the transpose of the matrix $U$.

Let $\cB_b(\mR^{d})$ be the set of all bounded Borel measurable functions on $\mR^d$. Let $C(\mR^d)$ be the set of all  functions which are continuous on $\mR^d$. Let $C_c(\mR^d)$ be the set of all continuous functions with compact support. 

For a locally compact Polish space $\mU$ endowed with a metric $\rho_\mU$, $\cM(\mU)$ denotes the family of all measures $\nu$ on $(\mU, \sB(\mU))$ such that $\nu(K)<\infty$ for every 
compact subset $K \subset \mU$. We endow $\mathcal{M}(\mU)$ with the weakest topology such that, for every $f \in C_c(\mU)$, the function $\mathcal{M}(\mU)\ni\nu \rightarrow \int_{\mU} f(u) \nu(\dif u)\in\mR$ is continuous. Thus, $\cM(\mU)$ 
is metrizable as a Polish space (see, e.g. \cite{bdv2}). Throughout the paper, we fix a $\nu \in \cM(\mU)$. 

Let $\mW=C([0,T],\mR^l)$ be the space of all continuous functions from $[0,T]$ to $\mR^l$, equipped with the uniform convergence topology. Let $D([0,T],\mR^l)$ be the space of right continuous functions with left limits from $[0,T]$ to $\mR^l$, equipped with the following metric: for any $X, Y \in D([0,T],\mR^l)$,
$$
\rho_{D([0,T],\mR^l)}(X, Y):=\inf _{\lambda \in {\bf \Lambda}}\left\{\sup _{0 \leq t \leq T}\left|X_t-Y_{\lambda(t)}\right|+\sup _{0 \leq s<t \leq T}\left|\log \frac{\lambda(t)-\lambda(s)}{t-s}\right|\right\},
$$
where
$$
{\bf \Lambda}=\left\{\lambda=\lambda(t): \begin{array}{l} \lambda \text { is strictly increasing, continuous on } t \in[0, T] \\ \text { such that } \lambda(0)=0, \lambda(T)=T \end{array}\right\}.
$$
From \cite[Chapter 15]{hwy}, we know that $\rho_{D([0,T],\mR^l)}$ induces the Skorohod topology on $D([0,T],\mR^l)$ and the space $\left(D([0,T],\mR^l), \rho_{D([0,T],\mR^l)}\right)$ is a Polish space.

Let $\mM=\cM([0,T]\times\mU)$ and $\mV=\mW\times\mM$. For any $\t>0$, let $\mP_\t$ be the unique probability measure on $(\mathbb{V}, \sB(\mathbb{V}))$ such that 

$(i)$ the canonical map $W: \mathbb{V} \mapsto \mathbb{W}, W(w, m):=w$ is a $l$-dimensional Brownian motion; 

$(ii)$ the canonical map $N: \mathbb{V} \mapsto \mathbb{M}, N(w, m):=m$ is a Poisson random measure with the intensity measure $\t Leb_T \otimes \nu$, where $Leb_T$ is the Lebesgue measures on $[0,T]$; 

$(iii)$ $\{W_t\}_{t \in[0, T]}$ and $\{N((0, t] \times B)-\t t\nu(B)\}_{t \in[0, T]}$ are $\sG_t$-martingales for $B \in \sB(\mU)$, where
$$
\sG_t:=\sigma\{\left(W_s, N((0, s] \times B)\right): s \in(0, t], B \in\sB(\mU)\}.
$$

In order to apply the weak convergence method to prove the uniform LDP, we introduce controlled random measures. Let $\bar\mM=\cM([0,T]\times\mU\times[0,\infty))$ and $\bar\mV=\mW\times\bar\mM$. Define $(\bar\mP, \{\bar\sG_t\})$ on $(\bar\mV,\sB(\bar\mV))$ analogous to $(\mP_\t,\{\sG_t\})$ by replacing $(N, \t Leb_T \otimes \nu)$ with $(\bar N, Leb_T \otimes \nu \otimes Leb_{\infty})$, where $Leb_{\infty}$ is the Lebesgue measures on $[0, \infty)$. Let $\left\{\bar{\sF}_t\right\}$ be the $\bar{\mathbb{P}}$-completion of $\left\{\bar\sG_t\right\}$ and $\bar{\mathcal{P}}$ the predictable $\sigma$-field on $[0, T] \times \bar\mV$ with respect to the filtration $\left\{\bar{\sF}_t\right\}_{t \in[0, T]}$. Set
\ce
&&\mathcal{A}_1:=\bigg\{\xi:[0, T] \times\bar\mV\mapsto \mathbb{R}^l \mid \xi \text { is } \bar{\mathcal{P}} \backslash \sB(\mathbb{R}^l) \text {-measurable,}~\text{and}~\int_0^T|\xi(s)|^2 \mathrm{~d} s<\infty,~\bar{\mathbb{P}} \text {-a.s.}\bigg\},\\
&&\mathcal{A}_2:=\left\{\zeta: [0,T]\times \bar\mV\times\mU\mapsto [0, \infty) \mid \zeta \text { is }(\bar{\mathcal{P}} \otimes \sB(\mU)) \backslash \sB([0, \infty)) \text {-measurable}\right\},\\
&&\mathcal{A}:=\cA_1\times\cA_2.
\de
Moreover, for $\zeta\in\mathcal{A}_2$, we define a controlled random measure $N^{\zeta}$ on $[0,T]\times\mU$ by
$$
N^{\zeta}((0, t] \times U)=\int_{(0, t]}\int_U\int_0^\infty 1_{[0, \zeta(s, u)]}(r)\bar N(\dif s\dif u\dif r), \quad t \in[0, T], \quad U \in \sB(\mU),
$$
with $\zeta$ selecting the intensity for the points at location $u$ and time $s$, in a possibly random but non-anticipating way. If $\zeta(s, \omega, u)=\t$ for any $s \in[0, T], \omega\in\bar\mV, u\in \mU$, we write $N^\t$ instead of $N^\zeta$.

\subsection{Maximal monotone operators}\label{mmo}

Given a multivalued operator $A: \mR^d\mapsto 2^{\mR^d}$, where $2^{\mR^d}$ stands for all the subsets of $\mR^d$. Set
$$
\cD(A):= \left\{x\in \mR^d: A(x) \ne \emptyset\right\}
$$
and
$$
Gr(A):= \left\{(x,y)\in \mR^{2d}:x \in \cD(A), ~ y\in A(x)\right\}.
$$

\bd\label{momope}
A multivalued operator $A$ is called monotone if 
$$
\langle x_1 - x_2, y_1 - y_2 \rangle \geq 0,\ \text{for\ any}\ (x_1,y_1), (x_2,y_2) \in Gr(A).
$$
A monotone operator $A$ is called maximal monotone if
$$
(x_1,y_1) \in Gr(A) \iff \langle x_1-x_2, y_1 -y_2 \rangle \geq 0, \forall (x_2,y_2) \in Gr(A).
$$
\ed

\bx\label{lscfun}
Assume that a function $\varphi:\mR^d\mapsto(-\infty, +\infty]$ is lower semicontinuous convex, and ${\rm Int}(Dom(\varphi))\neq \emptyset$, where $Dom(\varphi)\equiv\{x\in\mR^d; \varphi(x)<\infty\}$ and $\operatorname{Int}(Dom(\varphi))$ is the interior of $Dom(\varphi)$. Define the subdifferential operator of the function $\varphi$:
$$
\partial\varphi(x):=\{y\in\mR^d: \<y,z-x\>+\varphi(x)\leq \varphi(z), \forall z\in\mR^d\}.
$$
Then $\partial\varphi$ is a maximal monotone operator.
\ex

\bx\label{conset} 
Suppose that $\mathcal{O}$ is a closed convex subset of $\mathbb{R}^d$ and $\operatorname{Int}(\mathcal{O})\neq\emptyset$. Define
$$
\cI_{\mathcal{O}}(x):= \begin{cases}0, & \text { if } x \in \mathcal{O}, \\ 
+\infty, & \text { if } x \notin \mathcal{O}.\end{cases}
$$
The subdifferential operator of $\cI_{\mathcal{O}}$ is given by
$$
\begin{aligned}
\partial \cI_{\mathcal{O}}(x) & :=\left\{y \in \mathbb{R}^d:\langle y, x-z\rangle \geq 0, \forall z \in \mathcal{O}\right\} \\
& = \begin{cases}\emptyset, & \text { if } x \notin \mathcal{O}, \\
\{0\}, & \text { if } x \in \operatorname{Int}(\mathcal{O}), \\
E_x, & \text { if } x \in \partial \mathcal{O},\end{cases}
\end{aligned}
$$
where $E_x$ is the exterior normal cone at $x$. We know that $\partial \cI_{\mathcal{O}}$ is a maximal monotone operator.
\ex

Given $T>0$. Let $\sV_{0}$ be the set of all c\`adl\`ag functions $K: [0,T]\mapsto\mR^{d}$ with finite variations and $K_{0} = 0$. For $K\in\sV_0$ and $s\in [0,T]$, we shall use $|K|_{0}^{s}$ to denote the variation of $K$ on $[0,s]$. Set
\ce
&&\sA:=\Big\{(X,K): X\in D([0,T],\overline{\cD(A)}), K \in \sV_0, \\
&&\qquad\qquad\quad~\mbox{and}~\langle X_{t}-x, \dif K_{t}-y\dif t\rangle \geq 0 ~\mbox{for any}~ (x,y)\in Gr(A)\Big\}.
\de
We recall the following results.

\bl\label{equi}
For $X\in D([0,T],\overline{\cD(A)})$ and $K\in \sV_{0}$, the following statements are equivalent:
\begin{enumerate}[(i)]
\item $(X,K)\in \sA$.
\item For any $x,y\in D([0,T],\mR^d)$ with $(x_t,y_t)\in Gr(A)$, it holds that
$$
\left\langle X_t-x_t, \dif K_t-y_t\dif t\right\rangle \geq0.
$$
\item For any $(X^{'},K^{'})\in \sA$, it holds that
$$
\left\langle X_t-X_t^{'},\dif K_t-\dif K_t^{'}\right\rangle \geq0.
$$
\end{enumerate}
\el

\bl\label{inteineq}
Assume that the interior of the set $\cD(A)$ is non-empty, i.e. $\rm{Int}(\cD(A))\ne\emptyset$. For any $a\in \rm{Int}(\cD(A))$, there exist $M_1 >0$, and $M_{2},M_{3}\geq0$ such that  for any $(X,K)\in \sA$ and $0\leq s<t\leq T$,
$$
\int_s^t{\left< X_r-a, \dif K_r \right>}\geq M_1\left| K \right|_{s}^{t}-M_2\int_s^t{\left| X_r-a\right|}\dif r-M_3\left( t-s \right) .
$$
\el

\subsection{Uniform LDPs}

In this subsection, we introduce the uniform LDPs.

Let $(\mX, \rho_\mX)$ be a Polish space and let $\mX_0$ be a set. Let $\{X^{\e}(x_0), \e>0, x_0\in \mX_0\}$ be a family of $\mX$-valued random variables defined on $(\Omega, \mathscr{F}, \{\mathscr{F}_t\}_{t\in[0,T]}, \mP)$.

\bd\label{rfde} 
$(i)$ For $x_0\in\mX_0$, a function $\Lambda_{x_0}: \mX\rightarrow[0,+\infty]$ is called a rate function on $\mX$, if for all $M\geq 0$, $\{\phi\in\mX: \Lambda_{x_0}(\phi)\leq M\}$ is a closed subset of $\mX$.

$(ii)$ For $x_0\in\mX_0$, a function $\Lambda_{x_0}: \mX\rightarrow[0,+\infty]$ is called a good rate function on $\mX$, if for all $M\geq 0$, $\{\phi\in\mX: \Lambda_{x_0}(\phi)\leq M\}$ is a compact subset of $\mX$.
	
$(iii)$ A family $\{\Lambda_{x_0}, x_0\in\mX_0\}$ of rate functions is said to have compact level sets on compacts if for all compact subsets $B$ of $\mX_0$ and each $M\geq 0$, 
$\bigcup\limits_{x_0\in B}\left\{\phi \in \mX:\Lambda_{x_0}(\phi)\leq M\right\}$ is a compact subset of $\mX$.
\ed

\bd [Freidlin-Wentzell uniform LDP]\label{fwuldpde}
Let $\sA$ be a collection of subsets of $\mX_0$. $\{X^{\e}(x_0), \e>0, x_0\in \mX_0\}$ are said to satisfy a Freidlin-Wentzell uniform LDP with respect to the rate function $\Lambda_{x_0}$ with the speed $\e$ uniformly over $\sA$, if

$(a)$ For any $B\in \sA, M, \d, \t>0$, there exists a $\e_0>0$ such that
\ce
\mP(\rho_\mX(X^{\e}(x_0), \phi)<\d)\geq \exp\left\{-\frac{\Lambda_{x_0}(\phi)+\t}{\e}\right\},
\de
for all $0<\e<\e_0, x_0\in B, \phi\in \Phi_{x_0}(M)$, where $\Phi_{x_0}(M):=\{\phi\in\mX: \Lambda_{x_0}(\phi)\leq M \}$.

$(b)$ For any $B\in \sA, M, \d, \t>0$, there exists a $\e_0>0$ such that
\ce
\mP(\rho_\mX(X^{\e}(x_0), \Phi_{x_0}(M'))\geq\d)\leq \exp\left\{-\frac{M'-\t}{\e}\right\},
\de
for all $0<\e<\e_0, x_0\in B, 0\leq M'\leq M$, where 
\ce
\rho_\mX(X^{\e}(x_0), \Phi_{x_0}(M')):=\inf\limits_{\phi\in\Phi_{x_0}(M')} \rho_\mX(X^{\e}(x_0), \phi).
\de
\ed

\bd[Uniform Laplace principle]\label{ulpde} 
Let $\sA$ be a collection of subsets of $\mX_0$. $\{X^{\e}(x_0), \e>0, x_0\in \mX_0\}$ are said to satisfy a uniform Laplace principle with respect to the rate function $\Lambda_{x_0}$ with the speed $\e$ uniformly over $\sA$, if for any $B\in \sA$ and any bounded continuous $\Psi: \mX\rightarrow \mR$,
\ce
\lim\limits_{\e\rightarrow 0}\sup\limits_{x_0\in B}\left|\e \log \mE \exp\left(-\frac{\Psi(X^{\e}(x_0))}{\e}\right)+\inf\limits_{\phi\in \mX}\{\Psi(\phi)+\Lambda_{x_0}(\phi)\}\right|=0.
\de
\ed

If $\sA$ is the collection of compact subsets of $\mX_0$, and the family $\{\Lambda_{x_0}, x_0\in\mX_0\}$ of rate functions has compact level sets on compacts, Theorem 2.5 in \cite{ms} assures that  the Freidlin-Wentzell uniform LDP and the uniform Laplace principle are equivalent. Thus, in order to obtain that the family of $\{X^{\e}(x_0), \e>0, x_0\in \mX_0\}$ satisfies the Freidlin-Wentzell uniform LDP, we only need to prove that the family of $\{X^{\e}(x_0), \e>0, x_0\in \mX_0\}$ satisfies the uniform Laplace principle. Let us state the conditions under which the uniform Laplace principle holds. 

Set for $N \in \mathbb{N}$
$$
\mathbf{D}_1^N=\left\{h: [0, T]\rightarrow\mR^l| h ~\mbox{is}~ \sB([0,T])/\sB(\mR^l) ~\mbox{measurable, and}~ \int_{0}^{T}|h(t)|^{2}\dif t \leq N\right\},
$$
and we equip $\mathbf{D}_1^{N}$ with the weak convergence topology in $L^2\left([0, T], \mR^l\right)$. So, $\mathbf{D}_1^{N}$ is metrizable as a compact Polish space. Also, let
\ce
\mathbf{D}_2^N=\bigg\{g: [0,T]\times\mU\rightarrow[0, \infty)| g ~\mbox{is}~ \sB([0,T]\times\mU)/\sB([0, \infty)) ~\mbox{measurable, and}~ \\
 \int_0^T\int_{\mU}\ell\(g(t,u)\)\nu(\dif u)\dif t\leq N\bigg\},
\de
where $\ell:[0, \infty) \rightarrow[0, \infty)$ is defined by
$$
\ell(r)=r \log r-r+1, \quad r \in[0, \infty).
$$
Then we identify $g\in \mathbf{D}_2^N$ with a measure $\nu_T^g \in \cM([0,T]\times\mU)$, defined by 
$$
\nu_T^g(B):=\int_{B} g(s, u) \nu(\dif u)\dif s, \quad B \in \sB\left([0,T]\times\mU\right). 
$$
Throughout we consider the topology on $\mathbf{D}_2^N$ obtained through this identification, which makes $\mathbf{D}_2^N$ a compact space. Let $\mathbf{D}^N=\mathbf{D}_1^N\times\mathbf{D}_2^N$ with the usual product topology and $\mathbf{D}=\bigcup\limits_{N \geq 1}\mathbf{D}^N$. Let $\mathcal{A}^N$ be the space of $\mathbf{D}^N$-valued random controls:
$$
\cA^N=\left\{(\xi, \zeta) \in \mathcal{A}: (\xi(\omega), \zeta(\omega)) \in \mathbf{D}^N, \mathbb{P}\text{-a.s. } \omega\right\}.
$$

Let $\mathcal{G}^\e: \mX_0\times \mathbb{V}\rightarrow \mX$ be a measurable mapping for any $0<\e<1$.

\bco\label{cond}
There exists a measurable mapping $\mathcal{G}^0: \mX_0\times \mathbb{V}\rightarrow \mX$ such that the following holds.

$(i)$ Let $\left(\xi^\e, \zeta^\e\right), (\xi, \zeta) \in \mathcal{A}^N, x_0^\e, x_0\in \mX_0$ be such that, as $\e\rightarrow 0, \left(\xi^\e, \zeta^\e\right)$ converges in distribution to $(\xi, \zeta)$ and $x_0^\e\rightarrow x_0$. Then
$$
\mathcal{G}^\e\left(x_0^\e, \sqrt{\e} W+\int_0^{\cdot} \xi^\e(s) \dif s, \e N^{\e^{-1} \zeta^\e}\right) \Rightarrow \mathcal{G}^0\left(x_0, \int_0^{\cdot} \xi(s) \dif s, \nu_T^{\zeta}\right).
$$

$(ii)$ For $n \in \mathbb{N}$ let $\left(h_n, g_n\right), (h, g) \in \mathbf{D}^N, x_{0,n}, x_0\in\mX_0$ be such that $\left(x_{0,n}, h_n, g_n\right) \rightarrow(x_0, h, g)$ as $n\rightarrow \infty$. Then
$$
\mathcal{G}^0\left(x_{0,n}, \int_0^{\cdot} h_n(s)\dif s, \nu_T^{g_n}\right) \rightarrow \mathcal{G}^0\left(x_0, \int_0^{\cdot}h(s) \dif s, \nu_T^g\right).
$$
\eco

For any $x_0\in\mX_0$ and $\phi \in \mX$, set 
$$
\mathbf{D}_\phi:=\left\{(h, g) \in\mathbf{D}: \phi=\mathcal{G}^0\left(x_0, \int_0^{\cdot} h(s) \dif s, \nu_T^g\right)\right\}.
$$
Then we define
\ce
\Lambda_{x_0}(\phi)=\left\{\begin{array}{ll}
\inf\limits _{(h, g) \in\mathbf{D}_\phi}\left\{\frac{1}{2}\int_0^T|h(t)|^2\dif t+\int_0^T\int_{\mU}\ell\(g(t,u)\)\nu(\dif u)\dif t\right\}, \quad \mathbf{D}_\phi\neq\emptyset,\\
\infty, \qquad\qquad\qquad\qquad\qquad\qquad\qquad\qquad\qquad\qquad\qquad \mathbf{D}_\phi=\emptyset.
\end{array}
\right.
\de

The following result is due to \cite[Theorem 4.4]{Mar} and \cite[Theorem 2.5]{ms}.

\bt\label{ulpjuth}
Set $X^{\e}(x_0)=\mathcal{G}^\e\left(x_0, \sqrt{\e} W, \e N^{\e^{-1}}\right)$. Suppose that Condition \ref{cond} holds and for all $\phi \in \mX, x_0 \mapsto \Lambda_{x_0}(\phi)$ is a lower semicontinuous mapping from $\mX_0$ to $[0, \infty]$. Then, for all $x_0\in\mX_0, \phi\mapsto \Lambda_{x_0}(\phi)$ is a rate function on $\mX$ and the family $\{\Lambda_{x_0}, x_0\in\mX_0\}$ of rate functions has compact level sets on compacts. Furthermore, the family $\left\{X^{\e}(x_0)\right\}$ satisfies a Freidlin-Wentzell uniform LDP on $\mX$, with the rate function $\Lambda_{x_0}$, uniformly on compact subsets of $\mX_0$.
\et

\section{Main results}\label{main}

In this section, we formulate the main results in this paper. 

\subsection{Two uniform LDPs for multivalued SDEs with jumps}

In this subsection, we state the Freidlin-Wentzell uniform LDP and the Dembo-Zeitouni uniform LDP for Eq.(\ref{msdej1}).

We fix $T>0$ and recall Eq.(\ref{msdej1}), i.e.
\ce\left\{\begin{array}{ll}
\dif X^\e_t\in -A(X^\e_t)\dif t+b(X^\e_t)\dif t+\sqrt \e\sigma(X^\e_t)\dif W_t+\e\int_{\mU}f(X^\e_{t-},u) \tilde N^{\e^{-1}}(\dif t\dif u),\\
X^\e_0=x_0\in \overline{\cD(A)}, \quad 0\leq t\leq T.
\end{array}
\right.
\de

Assume:
\begin{enumerate}[$(\mathbf{H}_{A})$]
\item $0\in{\rm Int}(\cD(A))$ and $0\in A(0)$.
\end{enumerate}
\begin{enumerate}[$({\bf H}^{1}_{b,\s})$]
\item There exists a constant $L_1>0$ such that for $x,x'\in \overline{\cD(A)}$
$$
|b(x)-b(x')|+\|\sigma(x)-\sigma(x')\|\leq L_1|x-x'|.
$$
\end{enumerate}
\begin{enumerate}[$({\bf H}^{1}_{f})$]
\item For any $x\in \overline{\cD(A)}$, $u\in \mU$, $x+f(x,u)\in \overline{\cD(A)}$.
\end{enumerate}
\begin{enumerate}[$({\bf H}^{2}_{f})$]
\item There exists a positive function $L_2(u)$ satisfying
\ce
\sup\limits_{u\in\mU}L_2(u)\leq \gamma_1<1 ~\mbox{and}~ \int_{\mU}L^2_2(u)\nu(\dif u)<\infty,
\de
such that for $x,x'\in \overline{\cD(A)}$ and $u\in \mU$
\ce
|f(x,u)-f(x',u)|\leq L_2(u)|x-x'|,
\de
and
\ce
|f(0,u)|\leq L_2(u).
\de
\end{enumerate}
\begin{enumerate}[$({\bf H}^{3}_{f})$]
\item There exists a constant $\g_2>0$ such that for any $U\in \sB(\mU)$ with $\nu(U)<\infty$,
\ce
\int_{U}\exp\{\g_2 L^2_2(u)\}\nu(\dif u)<\infty.
\de
\end{enumerate}

\br
$(i)$ $(\mathbf{H}_{A})$ is equivalent to ${\rm Int}(\cD(A))\neq \emptyset$.

$(ii)$ $({\bf H}^{1}_{b,\s})$ implies that 
\ce
|b(x)|+\|\sigma(x)\|\leq (L_1+|b(0)|+\|\sigma(0)\|)(1+|x|).
\de

$(iii)$ $({\bf H}^{2}_{f})$ yields that
$$
\int_{\mU}|f(x,u)-f(x',u)|^2\nu(\dif u)\leq \int_{\mU}L^2_2(u)\nu(\dif u)|x-x'|^2,
$$
and
\ce
&&\int_{\mU}|f(x,u)|^2\nu(\dif u)\leq 2\int_{\mU}L^2_2(u)\nu(\dif u)(1+|x|^2), \\
&&\int_{\mU}|f(x,u)|^4\nu(\dif u)\leq 2^3\int_{\mU}L^2_2(u)\nu(\dif u)(1+|x|^4).
\de

$(iv)$ $(\mathbf{H}_{A})$, $({\bf H}^{1}_{b,\s})$, $({\bf H}^{1}_{f})$ and $({\bf H}^{2}_{f})$ assure the well-posedness for solutions of Eq.(\ref{msdej1}). And $({\bf H}^{3}_{f})$ is used to prove two uniform LDPs.
\er

Under $(\mathbf{H}_{A})$, $({\bf H}^{1}_{b,\s})$, $({\bf H}^{1}_{f})$ and $({\bf H}^{2}_{f})$, by Theorem 3.1 in \cite{gw}, we know that Eq.(\ref{msdej1}) has a unique strong solution $(X^\e(x_0), K^\e(x_0))$. That is, $(X^\e(x_0), K^\e(x_0))$ satisfies the following equation
\ce
X^\e_t(x_0)&=&x_0-K^\e_t(x_0)+\int_0^t b(X^\e_s(x_0))\dif s+\sqrt \e\int_0^t\sigma(X^\e_s(x_0))\dif W_s\\
&&+\e\int_0^t\int_{\mU}f(X^\e_{s-}(x_0),u) \tilde N^{\e^{-1}}(\dif s\dif u), \quad 0\leq t\leq T.
\de

Next, for $(h,g) \in\mathbf{D}$, consider the following equation:
\be\left\{\begin{array}{ll}
\dif X^{h, g}_t\in -A(X^{h, g}_t)\dif t+b(X^{h, g}_t)\dif t+\sigma(X^{h, g}_t)h(t)\dif t\\
\qquad\qquad+\int_{\mU}f(X^{h, g}_t,u)(g(t,u)-1)\nu(\dif u)\dif t, \quad 0\leq t\leq T,\\
X^{h, g}_0=x_0\in \overline{\cD(A)}.
\end{array}
\right.
\label{xx0eq}
\ee
Under $(\mathbf{H}_{A})$, $({\bf H}^{1}_{b,\s})$ and $({\bf H}^{2}_{f})$, by the penalization method, the above equation has a unique solution $(X^{h,g}(x_0), K^{h,g}(x_0))$. and set 
$$
\Lambda_{x_0}(\phi):=\inf _{(h,g) \in\mathbf{D}, \phi=X^{h,g}(x_0)}\left\{\frac{1}{2}\int_0^T|h(t)|^2\dif t+\int_0^T\int_{\mU}\ell\(g(t,u)\)\nu(\dif u)\dif t\right\}.
$$

Now, we state the Freidlin-Wentzell uniform LDP for Eq.(\ref{msdej1}), which is the first main result in this subsection.

\bt\label{fwuldpth}
Assume that $(\mathbf{H}_{A})$, $({\bf H}^{1}_{b,\s})$, $({\bf H}^{1}_{f})$-$({\bf H}^{3}_{f})$ hold. Then the family $\{X^\e(x_0), 0<\e<1, x_0\in\overline{\cD(A)}\}$ satisfies the Freidlin-Wentzell uniform LDP on $D([0,T],\overline{\cD(A)})$ with the good rate function $\Lambda_{x_0}$ uniformly on compact subsets of $\overline{\cD(A)}$.
\et

\br\label{fwthagst}
By the above theorem, we know that:

$(i)$ $\Lambda_{x_0}$ is a good rate function, that is, for all $M\geq 0$, $\Phi_{x_0}(M):=\{\phi\in D([0,T],\overline{\cD(A)}): \Lambda_{x_0}(\phi)\leq M\}$ is compact in $D([0,T],\overline{\cD(A)})$.

$(ii)$ For any compact subset $B$ of $\overline{\cD(A)}$, and any $M, \d, \t>0$, there exists a $\e_0>0$ such that
\ce
\mP(\rho_{D([0,T],\overline{\cD(A)})}(X^{\e}(x_0), \phi)<\d)\geq \exp\left\{-\frac{\Lambda_{x_0}(\phi)+\t}{\e}\right\},
\de
for all $0<\e<\e_0, x_0\in B, \phi\in \Phi_{x_0}(M)$.

$(iii)$ For any compact subset $B$ of $\overline{\cD(A)}$, and any $M, \d, \t>0$, there exists a $\e_0>0$ such that
\ce
\mP(\rho_{D([0,T],\overline{\cD(A)})}(X^{\e}(x_0), \Phi_{x_0}(M'))\geq\d)\leq \exp\left\{-\frac{M'-\t}{\e}\right\},
\de
for all $0<\e<\e_0, x_0\in B, 0\leq M'\leq M$.
\er

The proof of Theorem \ref{fwuldpth} is placed in Section \ref{uldpthproo}.

At present, we present the Dembo-Zeitouni uniform LDP for Eq.(\ref{msdej1}).

\bt\label{dzuldpth}
Assume that $(\mathbf{H}_{A})$, $({\bf H}^{1}_{b,\s})$, $({\bf H}^{1}_{f})$-$({\bf H}^{3}_{f})$ hold. Then the family $\{X^\e(x_0), 0<\e<1, x_0\in\overline{\cD(A)}\}$ satisfies the Dembo-Zeitouni uniform LDP on $D([0,T],\overline{\cD(A)})$ with the good rate function $\Lambda_{x_0}$ uniformly on compact subsets of $\overline{\cD(A)}$.
\et

\br
Theorem \ref{dzuldpth} implies that

$(a)$ For any compact subset $B$ of $\overline{\cD(A)}$ and any open $G \subset D([0,T],\overline{\cD(A)})$,
$$
\liminf _{\varepsilon \rightarrow 0} \inf _{x_0 \in B}\left(\e\log \mathbb{P}\left(X^{\varepsilon}(x_0) \in G\right)\right) \geq-\sup _{x_0 \in B} \inf_{\phi\in G}\Lambda_{x_0}(\phi).
$$

$(b)$ For any compact subset $B$ of $\overline{\cD(A)}$ and any closed $F \subset D([0,T],\overline{\cD(A)})$,
$$
\limsup _{\varepsilon \rightarrow 0} \sup _{x_0 \in B}\left(\e\log \mathbb{P}\left(X^{\varepsilon}(x_0) \in F\right)\right) \leq-\inf _{x_0 \in B} \inf_{\phi\in F}\Lambda_{x_0}(\phi).
$$
\er

The proof of Theorem \ref{dzuldpth} is postponed to Section \ref{uldpthproo}.

\subsection{The Freidlin-Wentzell LDP for invariant measures of multivalued SDEs with jumps}

In this subsection, we present the Freidlin-Wentzell LDP for invariant measures of Eq.(\ref{msdej1}).

Assume:
\begin{enumerate}[$({\bf H}^{2}_{\s})$]
\item There exists a constant $L_{\s}>0$ such that 
\ce
\sup\limits_{x\in\overline{\cD(A)}}\|\s(x)\|\leq L_{\s}.
\de
\end{enumerate}
\begin{enumerate}[$({\bf H}^{2'}_{f})$]
\item ``$|f(0,u)|\leq L_2(u)$" in $({\bf H}^{2}_{f})$ is replaced by ``$|f(x,u)|\leq L_2(u)$ for any $x\in\overline{\cD(A)}$", and other conditions in $({\bf H}^{2}_{f})$ remain the same.
\end{enumerate}
\begin{enumerate}[$({\bf H}_{b,\s,f})$]
\item There exists a constant $L_3>0$ such that for all $x\in\overline{\cD(A)}$ 
\ce
2\<x,b(x)\>+\|\sigma(x)\|^2+\int_{\mU}\big|f(x,u)\big|^2\nu(\dif u)\leq-L_3|x|^2. 
\de
\end{enumerate}

\br
$(i)$ $({\bf H}^{2'}_{f})$ is stronger than $({\bf H}^{2}_{f})$.

$(ii)$ $({\bf H}^{2}_{\s})$ and $({\bf H}^{2'}_{f})$ are used to obtain the exponential moment estimate for the solution of Eq.(\ref{msdej1}). And $({\bf H}_{b,\s,f})$ assures the existence of invariant probability measures for the solution of Eq.(\ref{msdej1}).
\er

Under $(\mathbf{H}_{A})$, $({\bf H}^{1}_{b,\s})$, $({\bf H}^{1}_{f})$, $({\bf H}^{2}_{f})$ and $({\bf H}_{b,\s,f})$, by \cite[Theorem 4.1]{gw}, we know that Eq.(\ref{msdej1}) has an invariant probability measure $\mu^\e$. Then we continue to study the LDP of $\{\mu^\e, 0<\e<1\}$. In order to do this, for any $T>0$ and $\phi\in D([0,T],\overline{\cD(A)})$, set
\ce
\Lambda_0^T(\phi):=\inf\left\{\frac{1}{2}\int_0^T|h(t)|^2\dif t+\int_0^T\int_{\mU}\ell\(g(t,u)\)\nu(\dif u)\dif t: (h, g) \in\mathbf{D}, \phi=X^{h,g}\right\},
\de
and it holds that $\Lambda_{x_0}=\Lambda_0^T$ with $\phi(0)=X^{h,g}_0(x_0)=x_0$. So, for $y\in\overline{\cD(A)}$, define the following function $V$
\be
V(y):=\inf\{\Lambda_0^T(\phi): \phi(0)=0, \phi(T)=y, T>0\},
\label{ratefuncinvi}
\ee
and it holds that the quasi-potential $V(y)\geq 0$ and $V(0)=0$ (\cite{fw}). The following conclusion describes the LDP estimate for $\{\mu^\e, 0<\e<1\}$.

\bt\label{ldpinvameas}
Suppose that $(\mathbf{H}_{A})$, $({\bf H}^{1}_{b,\s})$, $({\bf H}^{2}_{\s})$, $({\bf H}^{1}_{f})$, $({\bf H}^{2'}_{f})$, $({\bf H}^{3}_{f})$ and $({\bf H}_{b,\s,f})$ hold. Then $\{\mu^\e, 0<\e<1\}$ satisfies the Freidlin-Wentzell LDP on $D([0,T],\overline{\cD(A)})$ with the good rate function $V$.
\et

\br
Since the rate function $V$ is good, the Freidlin-Wentzell LDP for $\{\mu^\e, 0<\e<1\}$ is equivalent to the Dembo-Zeitouni LDP for $\{\mu^\e, 0<\e<1\}$. Thus, by Theorem \ref{ldpinvameas}, we know that $\{\mu^\e, 0<\e<1\}$ also satisfies the Dembo-Zeitouni LDP.
\er

The proof of the above theorem is postponed to Section \ref{ldpinvameasproo}.

\section{Proof of Theorem \ref{fwuldpth} and \ref{dzuldpth}}\label{uldpthproo}

In this section, we prove Theorem \ref{fwuldpth} and \ref{dzuldpth}. 

First of all, for the unique strong solution $(X^\e(x_0), K^\e(x_0))$ of Eq.(\ref{msdej1}), by the definition of strong solutions, there exists a unique measurable mapping $\mathcal{G}^\e: \overline{\cD(A)}\times \mathbb{V}\rightarrow D([0,T],\overline{\cD(A)})$ such that
$$
X^\e(x_0)=\mathcal{G}^\e\left(x_0, \sqrt{\e} W, \e N^{\e^{-1}}\right).
$$
Then for $\left(\xi^\e, \zeta^\e\right)\in \mathcal{A}^N$ and $x^\e_0\in \overline{\cD(A)}$, consider the following controlled equation:
\be\left\{\begin{array}{ll}
\dif X^{\e,\xi^\e, \zeta^\e}_t\in -A(X^{\e,\xi^\e, \zeta^\e}_t)\dif t+b(X^{\e,\xi^\e, \zeta^\e}_t)\dif t+\sigma(X^{\e,\xi^\e, \zeta^\e}_t)\xi^\e(t)\dif t\\
\qquad\qquad\qquad+\sqrt \e\sigma(X^{\e,\xi^\e, \zeta^\e}_t)\dif W_t+\int_{\mU}f(X^{\e,\xi^\e, \zeta^\e}_t,u)(\zeta^\e(t,u)-1)\nu(\dif u)\dif t\\
\qquad\qquad\qquad+\e\int_{\mU}f(X^{\e,\xi^\e, \zeta^\e}_{t-},u){\tilde N}^{\e^{-1}\zeta^\e}(\dif t\dif u),\\
X^{\e,\xi^\e, \zeta^\e}_0=x^\e_0\in \overline{\cD(A)}, \quad 0\leq t\leq T,
\end{array}
\right.
\label{skeleq}
\ee
where ${\tilde N}^{\e^{-1}\zeta^\e}(\dif t\dif u):=N^{\e^{-1}\zeta^\e}(\dif t\dif u)-\e^{-1}\zeta^\e(t,u)\nu(\dif u)\dif t$. Under $(\mathbf{H}_{A})$, $({\bf H}^{1}_{b,\s})$, $({\bf H}^{1}_{f})$ and $({\bf H}^{2}_{f})$, by Theorem 3.1 in \cite{gw}, we know that Eq.(\ref{skeleq}) has a unique strong solution $(X^{\e,\xi^\e, \zeta^\e}(x^\e_0), K^{\e,\xi^\e, \zeta^\e}(x^\e_0))$, that is,
\ce
X^{\e,\xi^\e, \zeta^\e}_t(x^\e_0)&=&x^\e_0-K^{\e,\xi^\e, \zeta^\e}_t(x^\e_0)+\int_0^tb(X^{\e,\xi^\e, \zeta^\e}_s(x^\e_0))\dif s+\int_0^t\sigma(X^{\e,\xi^\e, \zeta^\e}_s(x^\e_0))\xi^\e(s)\dif s\\
&&+\sqrt \e\int_0^t\sigma(X^{\e,\xi^\e, \zeta^\e}_s(x^\e_0))\dif W_s+\int_0^t\int_{\mU}f(X^{\e,\xi^\e, \zeta^\e}_s(x^\e_0),u)(\zeta^\e(s,u)-1)\nu(\dif u)\dif s\\
&&+\e\int_0^t\int_{\mU}f(X^{\e,\xi^\e, \zeta^\e}_{s-}(x^\e_0),u){\tilde N}^{\e^{-1}\zeta^\e}(\dif s\dif u).
\de
Moreover, by the Girsanov theorem for Brownian motions and random measures, it holds that 
\ce
X^{\e,\xi^\e, \zeta^\e}(x^\e_0)=\mathcal{G}^\e\left(x^\e_0, \sqrt{\e} W+\int_0^{\cdot} \xi^\e(s) \dif s, \e N^{\e^{-1} \zeta^\e}\right). 
\de

Next, for $\left(\xi, \zeta\right)\in \mathcal{A}^N$ and $x_0\in \overline{\cD(A)}$, we observe the following multivalued differential equation:
\be\left\{\begin{array}{ll}
\dif X^{\xi, \zeta}_t\in -A(X^{\xi, \zeta}_t)\dif t+b(X^{\xi, \zeta}_t)\dif t+\sigma(X^{\xi, \zeta}_t)\xi(t)\dif t\\
\qquad\qquad\qquad\qquad+\int_{\mU}f(X^{\xi, \zeta}_t,u)(\zeta(t,u)-1)\nu(\dif u)\dif t,\\
X^{\xi, \zeta}_0=x_0\in \overline{\cD(A)}, \quad 0\leq t\leq T.
\end{array}
\right.
\label{xx0eq2}
\ee
Under $(\mathbf{H}_{A})$, $({\bf H}^{1}_{b,\s})$ and $({\bf H}^{2}_{f})$, by the penalization method, the above equation has a unique solution $(X^{\xi, \zeta}(x_0),K^{\xi, \zeta}(x_0))$. That is, $(X^{\xi, \zeta}(x_0),K^{\xi, \zeta}(x_0))$ solves the following equation
\ce
X^{\xi, \zeta}_t(x_0)&=&x_0-K^{\xi, \zeta}_t(x_0)+\int_0^t b(X^{\xi, \zeta}_s(x_0))\dif s+\int_0^t \sigma(X^{\xi, \zeta}_s(x_0))\xi(s)\dif s\\
&&+\int_0^t\int_{\mU}f(X^{\xi, \zeta}_s(x_0),u)(\zeta(s,u)-1)\nu(\dif u)\dif s.
\de
So, we define a measurable mapping $\mathcal{G}^0: \overline{\cD(A)}\times \mathbb{V}\rightarrow D([0,T],\overline{\cD(A)})$ as follows
$$
X^{\xi, \zeta}(x_0)=\mathcal{G}^0\left(x_0, \int_0^{\cdot} \xi(s) \dif s, \nu_T^{\zeta}\right).
$$

In the following, we verify that $\cG^\e, \cG^0$ satisfy Condition \ref{cond} and prove Theorem \ref{fwuldpth} and \ref{dzuldpth}.

\subsection{Some estimates for $X^{\e,\xi^\e, \zeta^\e}, X^{\xi, \zeta}$}

In this subsection, we give some estimates for $X^{\e,\xi^\e, \zeta^\e}, X^{\xi, \zeta}$. 

First, let us prepare some results (cf. \cite[Remark 3.6]{bcd}). 

\bl
Suppose that $({\bf H}^{2}_{f})$ holds. Then it holds that
\be
&&\sup\limits_{\left(\xi, \zeta\right)\in \mathcal{A}^N}\int_0^T|\xi(s)|\dif s\leq (TN)^{1/2}, \quad \mP\mbox{-a.s.}, \label{xiboun}\\
&&\sup\limits_{\left(\xi, \zeta\right)\in \mathcal{A}^N}\int_0^T\int_{\mU}L^2_2(u)\zeta(s,u)\nu(\dif u)\dif s\leq M, \quad \mP\mbox{-a.s.},\label{l22zetaboun}\\
&&\sup\limits_{\left(\xi, \zeta\right)\in \mathcal{A}^N}\int_0^T\int_{\mU}L_2(u)|\zeta(s,u)-1|\nu(\dif u)\dif s\leq M, \quad \mP\mbox{-a.s.}, \label{l2zetaboun}
\ee
where the constant $M>0$ is non random, and
\be
\lim\limits_{\iota\rightarrow 0}\sup\limits_{\left(\xi, \zeta\right)\in \mathcal{A}^N}\sup\limits_{|t-s|\leq \iota}\int_s^t\int_{\mU}L_2(u)|\zeta(r,u)-1|\nu(\dif u)\dif r=0, \quad \mP\mbox{-a.s.}.
\label{l2zetalimi}
\ee
\el
\begin{proof}
First of all, by the H\"older inequality and $\xi(\omega)\in {\bf D}_1^N$ for $\mP$-a.s. $\omega$, it holds that
\ce
\int_0^T|\xi(s)|\dif s\leq\left(\int_0^T|\xi(s)|^2\dif s\right)^{1/2}T^{1/2}\leq (TN)^{1/2}.
\de
So, we get (\ref{xiboun}).

In the following, we treat $\int_0^T\int_{\mU}L^2_2(u)\zeta(s,u)\nu(\dif u)\dif s$. By the super linear growth of the function $\ell$, it holds that there exist two constants 
$\kappa_1, \kappa_2\in(0,\infty)$ such that for any $x\geq \kappa_1$, $x\leq \kappa_2 \ell(x)$. So,
\ce
\int_0^T\int_{\mU}L^2_2(u)\zeta(s,u)\nu(\dif u)\dif s&=&\int_{([0,T]\times\mU)\cap\{\zeta\geq \kappa_1\}}L^2_2(u)\zeta(s,u)\nu(\dif u)\dif s\\
&&+\int_{([0,T]\times\mU)\cap\{\zeta<\kappa_1\}}L^2_2(u)\zeta(s,u)\nu(\dif u)\dif s\\
&\leq&\int_{([0,T]\times\mU)\cap\{\zeta\geq \kappa_1\}}\kappa_2 \ell(\zeta(s,u))\nu(\dif u)\dif s\\
&&+\int_{([0,T]\times\mU)\cap\{\zeta<\kappa_1\}}L^2_2(u)\kappa_1\nu(\dif u)\dif s\\
&\leq&\kappa_2 N+\kappa_1T\int_\mU L^2_2(u)\nu(\dif u),
\de
where we use the assumption that
\ce
\sup\limits_{u\in\mU}L_2(u)\leq \gamma_1<1 ~\mbox{and}~ \int_{\mU}L^2_2(u)\nu(\dif u)<\infty.
\de
That is, (\ref{l22zetaboun}) is proved.

Next, we estimate $\int_s^t\int_{\mU}L_2(u)|\zeta(r,u)-1|\nu(\dif u)\dif r$ for $0\leq s<t\leq T$. Note that for any $\t>0$,
\ce
|x-1| \leq c_1(\t) \ell(x), \quad |x-1| >\t,\\
|x-1|^2 \leq c_2(\t) \ell(x), \quad |x-1| \leq \t,
\de
where $c_1(\t)>0$ and $c_1(\t) \rightarrow 0$ as $\t \rightarrow \infty$ and $0<c_2(\t)<\infty$ is a constant (cf. \cite[Remark 3.3]{bcd}). Thus, from the H\"older inequality, it follows that
\ce
&&\int_s^t\int_{\mU}L_2(u)|\zeta(r,u)-1|\nu(\dif u)\dif r\\
&=&\int_{([s,t]\times\mU)\cap\{|\zeta-1|>\t\}}L_2(u)|\zeta(r,u)-1|\nu(\dif u)\dif r+\int_{([s,t]\times\mU)\cap\{|\zeta-1|\leq\t\}}L_2(u)|\zeta(r,u)-1|\nu(\dif u)\dif r\\
&\leq&\int_{([s,t]\times\mU)\cap\{|\zeta-1|>\t\}}c_1(\t) \ell(\zeta(r,u))\nu(\dif u)\dif r\\
&&+\left(\int_{([s,t]\times\mU)\cap\{|\zeta-1|\leq\t\}}L^2_2(u)\nu(\dif u)\dif r\right)^{1/2}\left(\int_{([s,t]\times\mU)\cap\{|\zeta-1|\leq\t\}}|\zeta(r,u)-1|^2\nu(\dif u)\dif r\right)^{1/2}\\
&\leq&c_1(\t)N+\left(\int_{\mU}L^2_2(u)\nu(\dif u)\right)^{1/2}(t-s)^{1/2}c^{1/2}_2(\t)N^{1/2}.
\de

Taking $s=0, t=T$, we obtain that 
\ce
\int_0^T\int_{\mU}L_2(u)|\zeta(r,u)-1|\nu(\dif u)\dif r\leq c_1(\t)N+\left(\int_{\mU}L^2_2(u)\nu(\dif u)\right)^{1/2}T^{1/2}c^{1/2}_2(\t)N^{1/2}.
\de
This is just (\ref{l2zetaboun}).

Finally, we take $s,t\in[0,T], |t-s|<\iota$, and by the above deduction conclude (\ref{l2zetalimi}).
\end{proof}

In the following, we investigate $\{(X^{\e,\xi^\e, \zeta^\e}(x_0), K^{\e,\xi^\e, \zeta^\e}(x_0)), 0<\e<1\}$.

\bl\label{xekeboun}
Under the assumptions of Theorem \ref{fwuldpth}, it holds that 
\be
&&\sup\limits_{\e}\mE\sup\limits_{t\in[0,T]}|X^{\e,\xi^\e, \zeta^\e}_t(x_0)|^4\leq C(1+|x_0|^4),\label{xexezeboun}\\
&&\sup\limits_{\e}\mE|K^{\e,\xi^\e, \zeta^\e}(x_0)|_0^T\leq C(1+|x_0|^2).\label{kexezeboun}
\ee
\el
\begin{proof}
By the It\^o formula, it holds that for any $t\in[0,T]$
\be
&&|X^{\e,\xi^\e, \zeta^\e}_t(x_0)|^2\no\\
&=&|x_0|^2-2\int_0^t\<X^{\e,\xi^\e, \zeta^\e}_s(x_0), \dif K^{\e,\xi^\e, \zeta^\e}_s(x_0)\>+2\int_0^t\<X^{\e,\xi^\e, \zeta^\e}_s(x_0), b(X^{\e,\xi^\e, \zeta^\e}_s(x_0))\>\dif s\no\\
&&+2\int_0^t\<X^{\e,\xi^\e, \zeta^\e}_s(x_0), \sigma(X^{\e,\xi^\e, \zeta^\e}_s(x_0))\xi^\e(s)\>\dif s+2\sqrt \e\int_0^t\<X^{\e,\xi^\e, \zeta^\e}_s(x_0), \sigma(X^{\e,\xi^\e, \zeta^\e}_s(x_0))\dif W_s\>\no\\
&&+\e\int_0^t\|\sigma(X^{\e,\xi^\e, \zeta^\e}_s(x_0))\|^2\dif s+2\int_0^t\int_{\mU}\<X^{\e,\xi^\e, \zeta^\e}_s(x_0),f(X^{\e,\xi^\e, \zeta^\e}_s(x_0),u)(\zeta^\e(s,u)-1)\>\nu(\dif u)\dif s\no\\
&&+\int_0^t\int_{\mU}\[\left|X^{\e,\xi^\e, \zeta^\e}_{s-}(x_0)+\e f(X^{\e,\xi^\e, \zeta^\e}_{s-}(x_0),u)\right|^2-|X^{\e,\xi^\e, \zeta^\e}_{s-}(x_0)|^2\]{\tilde N}^{\e^{-1}\zeta^\e}(\dif s\dif u)\no\\
&&+\e\int_0^t\int_{\mU}\left|f(X^{\e,\xi^\e, \zeta^\e}_s(x_0),u)\right|^2\zeta^\e(s,u)\nu(\dif u)\dif s.
\label{itoxexize}
\ee
Then by Lemma \ref{equi} and $(\mathbf{H}_{A})$, it holds that 
\ce
-2\int_0^t\<X^{\e,\xi^\e, \zeta^\e}_s(x_0), \dif K^{\e,\xi^\e, \zeta^\e}_s(x_0)\>\leq 0,
\de
which together with $({\bf H}^{1}_{b,\s})$, $({\bf H}^{2}_{f})$ and the fact that $|x|\leq 1+|x|^2$ for any $x\in\mR^d$ implies that 
\ce
&&|X^{\e,\xi^\e, \zeta^\e}_t(x_0)|^2\\
&\leq&|x_0|^2+C\int_0^t(1+|X^{\e,\xi^\e, \zeta^\e}_s(x_0)|^2)\dif s+C\int_0^t|\xi^\e(s)|(1+|X^{\e,\xi^\e, \zeta^\e}_s(x_0)|^2)\dif s\\
&&+2\sqrt\e\left|\int_0^t\<X^{\e,\xi^\e, \zeta^\e}_s(x_0), \sigma(X^{\e,\xi^\e, \zeta^\e}_s(x_0))\dif W_s\>\right|\\
&&+2\int_0^t\left(\int_{\mU}L_2(u)|\zeta^\e(s,u)-1|\nu(\dif u)+\int_{\mU}L^2_2(u)\zeta^\e(s,u)\nu(\dif u)\right)(1+|X^{\e,\xi^\e, \zeta^\e}_s(x_0)|^2)\dif s\\
&&+\left|\int_0^t\int_{\mU}\[\left|X^{\e,\xi^\e, \zeta^\e}_{s-}(x_0)+\e f(X^{\e,\xi^\e, \zeta^\e}_{s-}(x_0),u)\right|^2-|X^{\e,\xi^\e, \zeta^\e}_{s-}(x_0)|^2\]{\tilde N}^{\e^{-1}\zeta^\e}(\dif s\dif u)\right|\\
&\leq&|x_0|^2+2\sqrt\e\left|\int_0^t\<X^{\e,\xi^\e, \zeta^\e}_s(x_0), \sigma(X^{\e,\xi^\e, \zeta^\e}_s(x_0))\dif W_s\>\right|\\
&&+\left|\int_0^t\int_{\mU}\[\left|X^{\e,\xi^\e, \zeta^\e}_{s-}(x_0)+\e f(X^{\e,\xi^\e, \zeta^\e}_{s-}(x_0),u)\right|^2-|X^{\e,\xi^\e, \zeta^\e}_{s-}(x_0)|^2\]{\tilde N}^{\e^{-1}\zeta^\e}(\dif s\dif u)\right|\\
&&+\int_0^t\bigg(C+C|\xi^\e(s)|+\left(\int_{\mU}L_2(u)|\zeta^\e(s,u)-1|\nu(\dif u)+\int_{\mU}L^2_2(u)\zeta^\e(s,u)\nu(\dif u)\right)\bigg)\\
&&\qquad\qquad \times(1+|X^{\e,\xi^\e, \zeta^\e}_s(x_0)|^2)\dif s.
\de
So, (\ref{l22zetaboun}), (\ref{l2zetaboun}) and the Gronwall inequality imply that
\ce
&&1+\sup\limits_{s\in[0,t]}|X^{\e,\xi^\e, \zeta^\e}_s(x_0)|^2\\
&\leq& C\bigg(1+|x_0|^2+2\sqrt\e\sup\limits_{s\in[0,T]}\left|\int_0^s\<X^{\e,\xi^\e, \zeta^\e}_r(x_0), \sigma(X^{\e,\xi^\e, \zeta^\e}_r(x_0))\dif W_r\>\right|\\
&&+\sup\limits_{s\in[0,T]}\left|\int_0^s\int_{\mU}\[\left|X^{\e,\xi^\e, \zeta^\e}_{r-}(x_0)+\e f(X^{\e,\xi^\e, \zeta^\e}_{r-}(x_0),u)\right|^2-|X^{\e,\xi^\e, \zeta^\e}_{r-}(x_0)|^2\]{\tilde N}^{\e^{-1}\zeta^\e}(\dif r\dif u)\right|\bigg)
\de

In the following, by the H\"older inequality and the BDG inequality, it holds that
\ce
&&1+\mE\sup\limits_{s\in[0,T]}|X^{\e,\xi^\e, \zeta^\e}_s(x_0)|^4\\
&\leq& C(1+|x_0|^4)+2C\e\mE\left(\sup\limits_{s\in[0,T]}\left|\int_0^s\<X^{\e,\xi^\e, \zeta^\e}_r(x_0), \sigma(X^{\e,\xi^\e, \zeta^\e}_r(x_0))\dif W_r\>\right|^2\right)\\
&&+C\mE\Bigg(\sup\limits_{s\in[0,T]}\bigg|\int_0^s\int_{\mU}\[\left|X^{\e,\xi^\e, \zeta^\e}_{r-}(x_0)+\e f(X^{\e,\xi^\e, \zeta^\e}_{r-}(x_0),u)\right|^2-|X^{\e,\xi^\e, \zeta^\e}_{r-}(x_0)|^2\]\\
&&\qquad\qquad\qquad\qquad\qquad{\tilde N}^{\e^{-1}\zeta^\e}(\dif r\dif u)\bigg|^2\Bigg)\\
&\leq& C(1+|x_0|^4)+C\e\mE\left(\int_0^T(1+|X^{\e,\xi^\e, \zeta^\e}_r(x_0)|^4)\dif r\right)\\
&&+C\e\mE\int_0^T\left(\int_{\mU}L^2_2(u)\zeta^\e(r,u)\nu(\dif u)\right)(1+|X^{\e,\xi^\e, \zeta^\e}_r(x_0)|^4)\dif r\\
&\leq& C(1+|x_0|^4)+\e C(T+M)(1+\mE\sup\limits_{r\in[0,T]}|X^{\e,\xi^\e, \zeta^\e}_r(x_0)|^4),
\de
which together with the arbitrariness of $\e$ yields that
\ce
\mE\sup\limits_{r\in[0,T]}|X^{\e,\xi^\e, \zeta^\e}_r(x_0)|^4\leq C(1+|x_0|^4).
\de
This is just (\ref{xexezeboun}).

Next, we show (\ref{kexezeboun}). From (\ref{itoxexize}) and Lemma \ref{inteineq}, it follows that
\ce
&&|X^{\e,\xi^\e, \zeta^\e}_T(x_0)|^2\no\\
&\leq&|x_0|^2-2M_1\left| K^{\e,\xi^\e, \zeta^\e}(x_0) \right|_{0}^{T}+M_2\int_0^T{\left|X^{\e,\xi^\e, \zeta^\e}_s(x_0)\right|}\dif s+M_3 T+C\int_0^T(1+|X^{\e,\xi^\e, \zeta^\e}_s(x_0)|^2)\dif s\\
&&+C\int_0^T|\xi^\e(s)|(1+|X^{\e,\xi^\e, \zeta^\e}_s(x_0)|^2)\dif s+2\sqrt \e\int_0^T\<X^{\e,\xi^\e, \zeta^\e}_s(x_0), \sigma(X^{\e,\xi^\e, \zeta^\e}_s(x_0))\dif W_s\>\\
&&+2\int_0^T\left(\int_{\mU}L_2(u)|\zeta^\e(s,u)-1|\nu(\dif u)+\int_{\mU}L^2_2(u)\zeta^\e(s,u)\nu(\dif u)\right)(1+|X^{\e,\xi^\e, \zeta^\e}_s(x_0)|^2)\dif s\\
&&+\int_0^T\int_{\mU}\[\left|X^{\e,\xi^\e, \zeta^\e}_{s-}(x_0)+\e f(X^{\e,\xi^\e, \zeta^\e}_{s-}(x_0),u)\right|^2-|X^{\e,\xi^\e, \zeta^\e}_{s-}(x_0)|^2\]{\tilde N}^{\e^{-1}\zeta^\e}(\dif s\dif u),
\de
and by taking the expectation on two sides
\ce
&&2M_1\mE\left| K^{\e,\xi^\e, \zeta^\e}(x_0) \right|_{0}^{T}\\
&\leq& |x_0|^2+M_3 T+C\mE\int_0^T(1+|X^{\e,\xi^\e, \zeta^\e}_s(x_0)|^2)\dif s+C\mE\int_0^T|\xi^\e(s)|(1+|X^{\e,\xi^\e, \zeta^\e}_s(x_0)|^2)\dif s\\
&&+2\mE\int_0^T\left(\int_{\mU}L_2(u)|\zeta^\e(s,u)-1|\nu(\dif u)+\int_{\mU}L^2_2(u)\zeta^\e(s,u)\nu(\dif u)\right)(1+|X^{\e,\xi^\e, \zeta^\e}_s(x_0)|^2)\dif s\\
&\leq& |x_0|^2+M_3 T+C(1+\mE\sup\limits_{s\in[0,T]}|X^{\e,\xi^\e, \zeta^\e}_s(x_0)|^2)\\
&\leq&C(1+ |x_0|^2).
\de
The proof is complete.
\end{proof}

Then the following proposition describes the property of $X^{\e,\xi^\e, \zeta^\e}(x)$ about the initial value $x$.

\bp\label{xex1x2es}
Under the assumptions of Theorem \ref{fwuldpth}, it holds that if for any $x^\e_0\in \overline{\cD(A)}$, $x_0^\e\rightarrow x_0$, then $X^{\e,\xi^\e, \zeta^\e}(x^\e_0)$ converges in probability to $X^{\e,\xi^\e, \zeta^\e}(x_0)$.
\ep
\begin{proof}
First of all, note that $X^{\e,\xi^\e, \zeta^\e}(x^\e_0), X^{\e,\xi^\e, \zeta^\e}(x_0)$ belong to $D([0,T],\overline{\cD(A)})$. Thus, we only need to prove that for any $\d>0$, and $\l^\e\in {\bf \Lambda}$ with $\lim\limits_{\e\rightarrow 0}\sup\limits_{t\in[0,T]}|\l^\e(t)-t|=0$ 
\ce
\lim\limits_{\e\rightarrow 0}\mP\left\{\sup\limits_{t\in[0,T]}|X^{\e,\xi^\e, \zeta^\e}_t(x^\e_0)-X^{\e,\xi^\e, \zeta^\e}_{\l^\e(t)}(x_0)|>\d\right\}=0.
\de
Note that
\ce
\mP\left\{\sup\limits_{t\in[0,T]}|X^{\e,\xi^\e, \zeta^\e}_t(x^\e_0)-X^{\e,\xi^\e, \zeta^\e}_{\l^\e(t)}(x_0)|>\d\right\}&\leq& \mP\left\{\sup\limits_{t\in[0,T]}|X^{\e,\xi^\e, \zeta^\e}_t(x^\e_0)-X^{\e,\xi^\e, \zeta^\e}_{t}(x_0)|>\d/2\right\}\\
&&+\mP\left\{\sup\limits_{t\in[0,T]}|X^{\e,\xi^\e, \zeta^\e}_t(x_0)-X^{\e,\xi^\e, \zeta^\e}_{\l^\e(t)}(x_0)|>\d/2\right\}.
\de
Thus, we divide into two steps to prove that the limits for two terms in the right side of the above inequality are zero as $\e$ tends to $0$.

{\bf Step 1.} We prove that 
\be
\lim\limits_{\e\rightarrow 0}\mP\left\{\sup\limits_{t\in[0,T]}|X^{\e,\xi^\e, \zeta^\e}_t(x^\e_0)-X^{\e,\xi^\e, \zeta^\e}_{t}(x_0)|>\d/2\right\}=0.
\label{xext}
\ee

Set $Z^{\e,\xi^\e, \zeta^\e}=X^{\e,\xi^\e, \zeta^\e}(x^\e_0)-X^{\e,\xi^\e, \zeta^\e}(x_0)$ and by the It\^o formula, Lemma \ref{equi}, the Taylor formula, $({\bf H}^{1}_{b,\s})$ and $({\bf H}^{2}_{f})$, we obtain that for any $t\in[0,T]$
\ce
&&|Z^{\e,\xi^\e, \zeta^\e}_t|^2\\
&=&|x^\e_0-x_0|^2-2\int_0^t\<Z^{\e,\xi^\e, \zeta^\e}_s, \dif (K^{\e,\xi^\e, \zeta^\e}_s(x^\e_0)-K^{\e,\xi^\e, \zeta^\e}_s(x_0))\>\\
&&+2\int_0^t\<Z^{\e,\xi^\e, \zeta^\e}_s, b(X^{\e,\xi^\e, \zeta^\e}_s(x^\e_0))-b(X^{\e,\xi^\e, \zeta^\e}_s(x_0))\>\dif s\\
&&+2\int_0^t\<Z^{\e,\xi^\e, \zeta^\e}_s, \[\sigma(X^{\e,\xi^\e, \zeta^\e}_s(x^\e_0))-\sigma(X^{\e,\xi^\e, \zeta^\e}_s(x_0))\]\xi^\e(s)\>\dif s\\
&&+2\int_0^t\int_{\mU}\<Z^{\e,\xi^\e, \zeta^\e}_s, \[f(X^{\e,\xi^\e, \zeta^\e}_s(x^\e_0),u)-f(X^{\e,\xi^\e, \zeta^\e}_s(x_0),u)\](\zeta^\e(s,u)-1)\>\nu(\dif u)\dif s\\
&&+\e\int_0^t\|\sigma(X^{\e,\xi^\e, \zeta^\e}_s(x^\e_0))-\sigma(X^{\e,\xi^\e, \zeta^\e}_s(x_0))\|^2\dif s\\
&&+\e\int_0^t\int_{\mU}\left|f(X^{\e,\xi^\e, \zeta^\e}_s(x^\e_0),u)-f(X^{\e,\xi^\e, \zeta^\e}_s(x_0),u)\)\right|^2\zeta^\e(s,u)\nu(\dif u)\dif s\\
&&+M_t^{\e,1}+M_t^{\e,2},
\de
where 
\ce
&&M_t^{\e,1}:=2\sqrt \e\int_0^t\<Z^{\e,\xi^\e, \zeta^\e}_s, \[\sigma(X^{\e,\xi^\e, \zeta^\e}_s(x^\e_0))-\sigma(X^{\e,\xi^\e, \zeta^\e}_s(x_0))\]\dif W_s\>,\\
&&M_t^{\e,2}:=\int_0^t\int_{\mU}\[\left|Z^{\e,\xi^\e, \zeta^\e}_{s-}+\e\(f(X^{\e,\xi^\e, \zeta^\e}_{s-}(x^\e_0),u)-f(X^{\e,\xi^\e, \zeta^\e}_{s-}(x_0),u)\)\right|^2-|Z^{\e,\xi^\e, \zeta^\e}_{s-}|^2\]\\
&&\qquad\qquad\qquad \times{\tilde N}^{\e^{-1}\zeta^\e}(\dif s\dif u).
\de
And Lemma \ref{equi} implies that
\ce
-2\int_0^t\<Z^{\e,\xi^\e, \zeta^\e}_s, \dif (K^{\e,\xi^\e, \zeta^\e}_s(x^\e_0)-K^{\e,\xi^\e, \zeta^\e}_s(x_0))\>\leq 0,
\de
which together with $({\bf H}^{1}_{b,\s})$ and $({\bf H}^{2}_{f})$ yields that
\ce
|Z^{\e,\xi^\e, \zeta^\e}_t|^2&\leq&|x^\e_0-x_0|^2+|M_t^{\e,1}|+|M_t^{\e,2}|+\int_0^t(2L_1+2L_1|\xi^\e(s)|+L_1^2)|Z^{\e,\xi^\e, \zeta^\e}_s|^2\dif s\\
&&+\int_0^t\left(\int_{\mU}2L_2(u)|\zeta^\e(s,u)-1|\nu(\dif u)+\int_{\mU}L^2_2(u)\zeta^\e(s,u)\nu(\dif u)\right)|Z^{\e,\xi^\e, \zeta^\e}_s|^2\dif s.
\de
Thus, (\ref{xiboun})-(\ref{l2zetaboun}) and the Gronwall inequality imply that
\ce
\sup\limits_{0\leq s\leq t}|Z^{\e,\xi^\e, \zeta^\e}_s|^2\leq C(|x^\e_0-x_0|^2+\sup\limits_{0\leq s\leq T}|M_s^{\e,1}|+\sup\limits_{0\leq s\leq T}|M_s^{\e,2}|).
\de

Next, from the BDG inequality, the mean value theorem, $({\bf H}^{1}_{b,\s})$, $({\bf H}^{2}_{f})$ and (\ref{l22zetaboun}), it follows that
\ce
\mE\sup\limits_{0\leq s\leq T}|Z^{\e,\xi^\e, \zeta^\e}_s|^4&\leq& C\left(|x^\e_0-x_0|^4+\mE\sup\limits_{0\leq s\leq T}|M_s^{\e,1}|^2+\mE\sup\limits_{0\leq s\leq T}|M_s^{\e,2}|^2\right)\\
&\leq&C|x^\e_0-x_0|^4+C\e\mE\int_0^T|Z^{\e,\xi^\e, \zeta^\e}_s|^4\dif s\\
&&+C\e\mE\int_0^T\left(\int_{\mU}L^2_2(u)\zeta^\e(s,u)\nu(\dif u)\right)|Z^{\e,\xi^\e, \zeta^\e}_s|^4\dif s\\
&\leq&C|x^\e_0-x_0|^4+\e(CT+CM)\mE\sup\limits_{0\leq s\leq T}|Z^{\e,\xi^\e, \zeta^\e}_s|^4.
\de
Since $\e$ is arbitrarily small, we conclude that
\be
\mE\sup\limits_{0\leq s\leq T}|Z^{\e,\xi^\e, \zeta^\e}_s|^4\leq C|x^\e_0-x_0|^4,
\label{zexiexeoxo}
\ee
which together with the Chebyshev inequality implies (\ref{xext}).

{\bf Step 2.} We prove that
\be
\lim\limits_{\e\rightarrow 0}\mP\left\{\sup\limits_{t\in[0,T]}|X^{\e,\xi^\e, \zeta^\e}_t(x_0)-X^{\e,\xi^\e, \zeta^\e}_{\l^\e(t)}(x_0)|>\d/2\right\}=0.
\label{xtlet}
\ee

Note that for $0\leq t<v<T$
\be
&&X^{\e,\xi^\e, \zeta^\e}_{v}(x_0)-X^{\e,\xi^\e, \zeta^\e}_{t}(x_0)\no\\
&=&-K^{\e,\xi^\e, \zeta^\e}_{v}(x_0)+K^{\e,\xi^\e, \zeta^\e}_{t}(x_0)+\int_t^v b(X^{\e,\xi^\e, \zeta^\e}_r(x_0))\dif r+\int_t^v\sigma(X^{\e,\xi^\e, \zeta^\e}_r(x_0))\xi^\e(r)\dif r\no\\
&&+\sqrt \e\int_t^v\sigma(X^{\e,\xi^\e, \zeta^\e}_r(x_0))\dif W_r+\int_t^v\int_{\mU}f(X^{\e,\xi^\e, \zeta^\e}_r(x_0),u)(\zeta^\e(r,u)-1)\nu(\dif u)\dif r\no\\
&&+\e\int_t^v\int_{\mU}f(X^{\e,\xi^\e, \zeta^\e}_{r-}(x_0),u){\tilde N}^{\e^{-1}\zeta^\e}(\dif r\dif u).
\label{xekeexpr}
\ee
Thus, by the It\^o formula, the H\"older inequality and $(\mathbf{H}^1_{b, \s})$ it holds that
\be
&&|X^{\e,\xi^\e, \zeta^\e}_{v}(x_0)-X^{\e,\xi^\e, \zeta^\e}_{t}(x_0)|^2\no\\
&=&-2\int_t^{v}\<X^{\e,\xi^\e, \zeta^\e}_{r}(x_0)-X^{\e,\xi^\e, \zeta^\e}_{t}(x_0), \dif K^{\e,\xi^\e, \zeta^\e}_{r}(x_0)\>\no\\
&&+2\int_t^{v}\<X^{\e,\xi^\e, \zeta^\e}_{r}(x_0)-X^{\e,\xi^\e, \zeta^\e}_{t}(x_0), b(X^{\e,\xi^\e, \zeta^\e}_{r}(x_0))\>\dif r\no\\
&&+2\int_t^{v}\<X^{\e,\xi^\e, \zeta^\e}_{r}(x_0)-X^{\e,\xi^\e, \zeta^\e}_{t}(x_0), \sigma(X^{\e,\xi^\e, \zeta^\e}_r(x_0))\xi^\e(r)\>\dif r\no\\
&&+2\sqrt \e\int_t^{v}\<X^{\e,\xi^\e, \zeta^\e}_{r}(x_0)-X^{\e,\xi^\e, \zeta^\e}_{t}(x_0), \sigma(X^{\e,\xi^\e, \zeta^\e}_r(x_0))\dif W_r\>+\e\int_t^{v}\|\sigma(X^{\e,\xi^\e, \zeta^\e}_r(x_0))\|^2\dif r\no\\
&&+2\int_t^{v}\int_{\mU}\<X^{\e,\xi^\e, \zeta^\e}_{r}(x_0)-X^{\e,\xi^\e, \zeta^\e}_{t}(x_0), f(X^{\e,\xi^\e, \zeta^\e}_r(x_0),u)(\zeta^\e(r,u)-1)\>\nu(\dif u)\dif r\no\\
&&+\int_t^{v}\int_{\mU}[|X^{\e,\xi^\e, \zeta^\e}_{r}(x_0)-X^{\e,\xi^\e, \zeta^\e}_{t}(x_0)+\e f(X^{\e,\xi^\e, \zeta^\e}_{r-}(x_0),u)|^2\no\\
&&\qquad\qquad -|X^{\e,\xi^\e, \zeta^\e}_{r}(x_0)-X^{\e,\xi^\e, \zeta^\e}_{t}(x_0)|^2]{\tilde N}^{\e^{-1}\zeta^\e}(\dif r\dif u)\no\\
&&+\e\int_t^{v}\int_{\mU}|f(X^{\e,\xi^\e, \zeta^\e}_{r-}(x_0),u)|^2\zeta^\e(r,u)\nu(\dif u)\dif r\no\\
&\leq&-2\int_t^{v}\<X^{\e,\xi^\e, \zeta^\e}_{r}(x_0)-X^{\e,\xi^\e, \zeta^\e}_{t}(x_0), \dif K^{\e,\xi^\e, \zeta^\e}_{r}(x_0)\>+\int_t^{v}|X^{\e,\xi^\e, \zeta^\e}_{r}(x_0)-X^{\e,\xi^\e, \zeta^\e}_{t}(x_0)|^2\dif r\no\\
&&+C\int_t^{v}(1+|X_{r}^{\e,\xi^\e,\zeta^\e}(x_0)|^2)\dif r+2\int_t^{v}|X^{\e,\xi^\e, \zeta^\e}_{r}(x_0)-X^{\e,\xi^\e, \zeta^\e}_{t}(x_0)|\|\sigma(X^{\e,\xi^\e, \zeta^\e}_r(x_0))\||\xi^\e(r)|\dif r\no\\
&&+2\int_t^{v}|X^{\e,\xi^\e, \zeta^\e}_{r}(x_0)-X^{\e,\xi^\e, \zeta^\e}_{t}(x_0)|\int_{\mU}|f(X^{\e,\xi^\e, \zeta^\e}_r(x_0),u)||\zeta^\e(r,u)-1|\nu(\dif u)\dif r\no\\
&&+2\sqrt \e\left|\int_t^{v}\<X^{\e,\xi^\e, \zeta^\e}_{r}(x_0)-X^{\e,\xi^\e, \zeta^\e}_{t}(x_0), \sigma(X^{\e,\xi^\e, \zeta^\e}_r(x_0))\dif W_r\>\right|\no\\
&&+\bigg{|}\int_t^{v}\int_{\mU}[|X^{\e,\xi^\e, \zeta^\e}_{r}(x_0)-X^{\e,\xi^\e, \zeta^\e}_{t}(x_0)+\e f(X^{\e,\xi^\e, \zeta^\e}_{r-}(x_0),u)|^2\no\\
&&\qquad\qquad -|X^{\e,\xi^\e, \zeta^\e}_{r}(x_0)-X^{\e,\xi^\e, \zeta^\e}_{t}(x_0)|^2]{\tilde N}^{\e^{-1}\zeta^\e}(\dif r\dif u)\bigg{|}\no\\
&&+\e(1+\sup\limits_{r\in[0,T]}|X_{r}^{\e,\xi^\e,\zeta^\e}|^2)M.
\label{xeuxes}
\ee

Next, we compute $-2\int_t^{v}\<X^{\e,\xi^\e, \zeta^\e}_{r}(x_0)-X^{\e,\xi^\e, \zeta^\e}_{t}(x_0), \dif K^{\e,\xi^\e, \zeta^\e}_{r}(x_0)\>$. Note that $0\in{\rm Int}(\cD(A))$. Thus, there is a $\t_0>0$ such that for any $R>0$ and $\t<\t_0$,
$$
\left\{y \in B(R): \rho_{\mR^d}(y,(\overline{\cD(A)})^c) \geqslant \t\right\} \neq \emptyset,
$$
where $B(R):=\{y\in\mR^d: |y|\leq R\}$ and $(\overline{\cD(A)})^c$ is the complement of $\overline{\cD(A)}$. Set
$$
\rho_R(\t):=\sup \left\{|z|: z \in A(y) \text { for all } y \in B(R) \text { with } \rho_{\mR^d}\left(y,(\overline{\cD(A)})^c\right) \geqslant \t\right\},
$$
and by the local boundedness of $A$ on ${\rm Int}(\cD(A))$ it holds that
$$
\rho_R(\t)<+\infty.
$$

Set for any $\iota>0$
$$
\vartheta_R(\iota):=\inf \left\{\t \in\left(0, \t_0\right): \rho_R(\t) \leqslant \iota^{-1 / 2}\right\}, 
$$
and it holds that
$$
\rho_R\left(\iota+\vartheta_R(\iota)\right) \leqslant \iota^{-1 / 2} \text { and } \quad \lim _{\iota \downarrow 0} \vartheta_R(\iota)=0.
$$
Take $\iota_R>0$ such that $\iota_R+\vartheta_R\left(\iota_R\right)<\t_0$. For $0<\iota<\iota_R \wedge 1$, let $X_t^{\e, \xi^\e, \zeta^\e, \iota, R}(x_0)$ be the projection of $X_t^{\e,\xi^\e,\zeta^\e}(x_0)$ on $\left\{y \in B(R): \rho_{\mR^d}\left(y,(\overline{\cD(A)})^c\right) \geqslant \iota+\vartheta_R(\iota)\right\}$. Thus, for $Y_t^{\e, \xi^\e, \zeta^\e, \iota, R} \in A(X_t^{\e, \xi^\e, \zeta^\e, \iota, R}(x_0)), \sup\limits_{v\in[0,T]}|X_v^{\e,\xi^\e,\zeta^\e}(x_0)| \leqslant R$ and $0<v-t<\iota$, it holds that
\ce
&&-2 \int_t^v\left\langle X_r^{\e,\xi^\e,\zeta^\e}(x_0)-X_t^{\e,\xi^\e,\zeta^\e}(x_0), \dif K_r^{\e,\xi^\e,\zeta^\e}(x_0) \right\rangle\\
 &=&-2 \int_t^v\left\langle X_r^{\e,\xi^\e,\zeta^\e}(x_0)-X_t^{\e, \xi^\e, \zeta^\e, \iota, R}(x_0), \dif K_r^{\e,\xi^\e,\zeta^\e}(x_0) \right\rangle\\
 &&-2 \int_t^v\left\langle X_t^{\e, \xi^\e, \zeta^\e, \iota, R}(x_0)-X_t^{\e,\xi^\e,\zeta^\e}(x_0), \dif K_r^{\e,\xi^\e,\zeta^\e}(x_0) \right\rangle \\ 
&\leqslant& -2 \int_t^v\left\langle X_r^{\e,\xi^\e,\zeta^\e}(x_0)-X_t^{\e, \xi^\e, \zeta^\e, \iota, R}(x_0), Y_t^{\e, \xi^\e, \zeta^\e, \iota, R}\right\rangle \dif r\\
&&+2\left(\iota+\vartheta_R(\iota)\right)\left|K^{\e,\xi^\e,\zeta^\e}(x_0) \right|_0^T\\
&\leqslant& 4 \iota^{1 / 2} R+2\left(\iota+\vartheta_R(\iota)\right)\left|K^{\e,\xi^\e,\zeta}(x_0) \right|_0^T,
\de
and furthermore by (\ref{xeuxes})
\be
&&\sup _{t \leqslant v\leqslant t+\iota}\left|X_v^{\e,\xi^\e,\zeta^\e}(x_0)-X_t^{\e,\xi^\e,\zeta^\e}(x_0)\right|^2 I_{\{\sup\limits_{v\in[0,T]}|X_v^{\e,\xi^\e,\zeta^\e}(x_0)| \leqslant R\}} \no\\
&\leqslant&\left(4 \iota^{1 / 2} R+2\left(\iota+\vartheta_R(\iota)\right)\left|K^{\e,\xi^\e,\zeta}(x_0) \right|_0^T\right)+4R^2\iota+C(1+R^2)\iota+\e(1+R^2)M\no\\
&&+C(1+R^2)\int_t^{t+\iota}|\xi^\e(r)|\dif r+8(1+R^2)\int_t^{t+\iota}\int_{\mU}L_2(u)|\zeta^\e(r,u)-1|\nu(\dif u)\dif r\no\\
&&+2\sqrt\e R\sup_{t \leqslant v\leqslant t+\iota}\left|\int_t^v\sigma(X^{\e,\xi^\e, \zeta^\e}_r(x_0))\dif W_r\right|\no\\
&&+2\sqrt\e\sup_{t \leqslant v\leqslant t+\iota}\left|\int_t^v\<X^{\e,\xi^\e, \zeta^\e}_r(x_0),\sigma(X^{\e,\xi^\e, \zeta^\e}_r(x_0))\dif W_r\>\right|\no\\
&&+2\e R\sup_{t \leqslant v\leqslant t+\iota}\left|\int_t^{v}\int_{\mU}f(X^{\e,\xi^\e, \zeta^\e}_{r-}(x_0),u){\tilde N}^{\e^{-1}\zeta^\e}(\dif r\dif u)\right|\no\\
&&+2\e \sup_{t \leqslant v\leqslant t+\iota}\left|\int_t^{v}\int_{\mU}\<X^{\e,\xi^\e, \zeta^\e}_{r-}(x_0), f(X^{\e,\xi^\e, \zeta^\e}_{r-}(x_0),u)\>{\tilde N}^{\e^{-1}\zeta^\e}(\dif r\dif u)\right|\no\\
&&+2\e\sup_{t \leqslant v\leqslant t+\iota}\left|\int_t^{v}\int_{\mU}|f(X^{\e,\xi^\e, \zeta^\e}_{r-}(x_0),u)|^2{\tilde N}^{\e^{-1}\zeta^\e}(\dif r\dif u)\right|.
\label{xres}
\ee

Besides, note that $\lim\limits_{\e\rightarrow 0}\sup\limits_{t\in[0,T]}|\l^\e(t)-t|=0$. Thus, there exists a $\e_0>0$ such that for any $\e<\e_0$ and $t\in[0,T]$, $|\l^\e(t)-t|<\iota$. And by the Chebyshev inequality it holds that
\be
&&\mP\left\{\sup\limits_{t\in[0,T]}|X^{\e,\xi^\e, \zeta^\e}_t(x_0)-X^{\e,\xi^\e, \zeta^\e}_{\l^\e(t)}(x_0)|>\d/2\right\}\no\\
&=&\mP\left\{\sup\limits_{t\in[0,T]}|X^{\e,\xi^\e, \zeta^\e}_t(x_0)-X^{\e,\xi^\e, \zeta^\e}_{\l^\e(t)}(x_0)|>\d/2,\sup\limits_{v\in[0,T]}|X_v^{\e,\xi^\e,\zeta^\e}(x_0)| \leqslant R\right\}\no\\
&&+\mP\left\{\sup\limits_{t\in[0,T]}|X^{\e,\xi^\e, \zeta^\e}_t(x_0)-X^{\e,\xi^\e, \zeta^\e}_{\l^\e(t)}(x_0)|>\d/2,\sup\limits_{v\in[0,T]}|X_v^{\e,\xi^\e,\zeta^\e}(x_0)|>R\right\}\no\\
&\leq&\mP\left\{\sup _{t\leqslant v \leqslant t+|\l^\e(t)-t|}\left|X^{\e,\xi^\e, \zeta^\e}_v(x_0)-X^{\e,\xi^\e, \zeta^\e}_t(x_0)\right|^2 I_{\{\sup\limits_{v\in[0,T]}|X_v^{\e,\xi^\e,\zeta^\e}(x_0)| \leqslant R\}}>\d^2/4\right\}\no\\
&&+\frac{1}{R^2}\mE\sup\limits_{v\in[0,T]}|X_v^{\e,\xi^\e,\zeta^\e}(x_0)|^2.
\label{mpxtlet}
\ee
In terms of (\ref{xres}), we estimate the first term of the right side for the above inequality. By (\ref{kexezeboun}) and the Chebyshev inequality, it holds that 
\ce
&&\mP\bigg\{C(1+R^2)|\l^\e(t)-t|^{1 / 2}+2\left(|\l^\e(t)-t|+\vartheta_R(|\l^\e(t)-t|)\right)\left|K^{\e,\xi^\e,\zeta}(x_0) \right|_0^T+\e(1+R^2)M\\
&&\qquad\quad >\d^2/28\bigg\}\\
&\leq&\frac{28}{\d^2}\left[C(1+R^2)|\l^\e(t)-t|^{1 / 2}+2C\left(|\l^\e(t)-t|+\vartheta_R(|\l^\e(t)-t|)\right)(1+|x_0|^2)+\e(1+R^2)M\right].
\de
And the integrability of $\xi^\e$ and (\ref{l2zetalimi}) yield that
\ce
&&\lim\limits_{\e\rightarrow 0}\mP\bigg\{C(1+R^2)\int_t^{t+|\l^\e(t)-t|}|\xi(r)|\dif r+8(1+R^2)\int_t^{t+|\l^\e(t)-t|}\int_{\mU}L_2(u)|\zeta^\e(r,u)-1|\nu(\dif u)\dif r\\
&&\qquad\quad >\d^2/28\bigg\}\\
&=&0.
\de

Then from the Burkholder-Davis-Gundy inequality and (\ref{xexezeboun}), it follows that
\ce
&&\mP\bigg\{2\sqrt\e R\sup_{t \leqslant v\leqslant t+|\l^\e(t)-t|}\left|\int_t^v\sigma(X^{\e,\xi^\e, \zeta^\e}_r(x_0))\dif W_r\right|>\d^2/28\bigg\}\\
&\leq&\frac{C}{\d^2}\sqrt\e R\left(\int_t^{t+|\l^\e(t)-t|}\mE\|\sigma(X^{\e,\xi^\e, \zeta^\e}_r(x_0))\|^2\dif r\right)^{1/2}\\
&\leq&\frac{C}{\d^2}\sqrt\e R(1+|x_0|)|\l^\e(t)-t|^{1/2},
\de
and furthermore
\ce
\lim\limits_{\e\rightarrow 0}\mP\bigg\{2\sqrt\e R\sup_{t \leqslant v\leqslant t+|\l^\e(t)-t|}\left|\int_t^v\sigma(X^{\e,\xi^\e, \zeta^\e}_r(x_0))\dif W_r\right|>\d^2/28\bigg\}=0.
\de
By the same deduction to that for the above limit, one can obtain that
\ce
\lim\limits_{\e\rightarrow 0}\mP\bigg\{2\sqrt\e\sup_{t \leqslant v\leqslant t+|\l^\e(t)-t|}\left|\int_t^v\<X^{\e,\xi^\e, \zeta^\e}_r(x_0),\sigma(X^{\e,\xi^\e, \zeta^\e}_r(x_0))\dif W_r\>\right|>\d^2/28\bigg\}=0.
\de

By the Chebyshev inequality and the Burkholder-Davis-Gundy inequality, it holds that
\ce
&&\mP\bigg\{2\e R\sup_{t \leqslant v\leqslant t+|\l^\e(t)-t|}\left|\int_t^{v}\int_{\mU}f(X^{\e,\xi^\e, \zeta^\e}_{r-}(x_0),u){\tilde N}^{\e^{-1}\zeta^\e}(\dif r\dif u)\right|>\d^2/28\bigg\}\\
&\leq&\frac{56^2}{\d^4}\e^2 R^2\mE\int_t^{t+|\l^\e(t)-t|}\int_{\mU}|f(X^{\e,\xi^\e, \zeta^\e}_{r-}(x_0),u)|^2\e^{-1}\zeta^\e(r,u)\nu(\dif u)\dif r\\
&\leq&\frac{2\times 56^2}{\d^4}\e R^2\mE(1+\sup\limits_{r\in[0,T]}|X_{r}^{\e,\xi^\e,\zeta^\e}|^2)M,
\de
and
\ce
\lim\limits_{\e\rightarrow 0}\mP\bigg\{2\e R\sup_{t \leqslant v\leqslant t+|\l^\e(t)-t|}\left|\int_t^{v}\int_{\mU}f(X^{\e,\xi^\e, \zeta^\e}_{r-}(x_0),u){\tilde N}^{\e^{-1}\zeta^\e}(\dif r\dif u)\right|>\d^2/28\bigg\}=0.
\de
Following the similar line to that for the above limit, we conclude that 
\ce
&&\lim\limits_{\e\rightarrow 0}\mP\bigg\{2\e \sup_{t \leqslant v\leqslant t+|\l^\e(t)-t|}\left|\int_t^{v}\int_{\mU}\<X^{\e,\xi^\e, \zeta^\e}_{r-}(x_0), f(X^{\e,\xi^\e, \zeta^\e}_{r-}(x_0),u)\>{\tilde N}^{\e^{-1}\zeta^\e}(\dif r\dif u)\right|>\d^2/28\bigg\}=0,\\
&&\lim\limits_{\e\rightarrow 0}\mP\bigg\{2\e\sup_{t \leqslant v\leqslant t+|\l^\e(t)-t|}\left|\int_t^{v}\int_{\mU}|f(X^{\e,\xi^\e, \zeta^\e}_{r-}(x_0),u)|^2{\tilde N}^{\e^{-1}\zeta^\e}(\dif r\dif u)\right|>\d^2/28\bigg\}=0.
\de

Finally, combining the above deduction with (\ref{mpxtlet}) and (\ref{xres}), and letting $\e\rightarrow 0$ first and then $R\rightarrow\infty$, one can obtain (\ref{xtlet}), which completes the proof.
\end{proof}

By the same deduction to that for (\ref{xexezeboun}) and (\ref{kexezeboun}), we can obtain the following result.

\bl\label{xkboun}
Under the assumptions of Theorem \ref{fwuldpth}, it holds that for $\left(\xi, \zeta\right)\in \mathcal{A}^N$
\be
&&\sup\limits_{t\in[0,T]}|X^{\xi, \zeta}_t(x_0)|^2\leq C(1+|x_0|^2), \mP\mbox{-a.s.}, \label{xxizeboun}\\
&&|K^{\xi, \zeta}(x_0)|_0^T\leq C(1+|x_0|^2), \mP\mbox{-a.s.},\label{kxizeboun}
\ee
where the constant $C>0$ is non random.
\el

\bl\label{zeezecon}
Suppose that the assumptions of Theorem \ref{fwuldpth} hold. Assume that $\left(\xi^\e, \zeta^\e\right), (\xi, \zeta) \in \mathcal{A}^N$, and $\e\rightarrow 0$, $\left(\xi^\e, \zeta^\e\right)$ converges almost surely to $(\xi, \zeta)$. Set
\ce
&&\Sigma_\e(t):=\int_0^t \s(X^{\xi,\zeta}_r(x_0))(\xi^\e(r)-\xi(r))\dif r,\\
&&\Psi_\e(t):=\int_0^t\int_{\mU}f(X^{\xi,\zeta}_r(x_0),u)\((\zeta^\e(r,u)-1)-(\zeta(r,u)-1)\)\nu(\dif u)\dif r.
\de
Then it holds that
\be
&&\lim\limits_{\e\rightarrow0}\sup\limits_{t\in[0,T]}|\Sigma_\e(t)|=0, \quad \mP\mbox{-a.s.},\label{gelimi}\\
&&\lim\limits_{\e\rightarrow 0}\sup\limits_{t\in[0,T]}|\Psi_\e(t)|=0, \quad \mP\mbox{-a.s.}. \label{psilimi}
\ee
\el
\begin{proof}
First of all, by the similar deduction to that in \cite[Lemma 6.4]{q1}, we can obtain (\ref{gelimi}). 

Next, we prove (\ref{psilimi}). By $({\bf H}^{2}_{f})$, (\ref{xxizeboun}) and (\ref{l2zetaboun}), it holds that
\ce
\sup\limits_{t\in[0,T]}|\Psi_\e(t)|&\leq& \int_0^T\int_{\mU}|f(X^{\xi,\zeta}_r(x_0),u)|\(|\zeta^\e(r,u)-1|+|\zeta(r,u)-1|\)\nu(\dif u)\dif r\\
&\leq&(1+\sup\limits_{t\in[0,T]}|X^{\xi, \zeta}_t(x_0)|)\int_0^T\int_{\mU}L_2(u)\(|\zeta^\e(r,u)-1|+|\zeta(r,u)-1|\)\nu(\dif u)\dif r\\
&\leq&C(1+|x_0|),
\de
where the constant $C>0$ is independent of $\e$. Then for $0\leq s<t\leq T$ and $|t-s|<\d$, $({\bf H}^{2}_{f})$ and (\ref{xxizeboun}) imply that 
\ce
|\Psi_\e(t)-\Psi_\e(s)|&\leq& \int_s^t\int_{\mU}|f(X^{\xi,\zeta}_r(x_0),u)|\(|\zeta^\e(r,u)-1|+|\zeta(r,u)-1|\)\nu(\dif u)\dif r\\
&\leq& (1+\sup\limits_{t\in[0,T]}|X^{\xi, \zeta}_t(x_0)|)\int_s^t\int_{\mU}L_2(u)\(|\zeta^\e(r,u)-1|+|\zeta(r,u)-1|\)\nu(\dif u)\dif r\\
&\leq&C(1+|x_0|)\bigg(\int_s^t\int_{\mU}L_2(u)|\zeta^\e(r,u)-1|\nu(\dif u)\dif r\\
&&\qquad\qquad\qquad +\int_s^t\int_{\mU}L_2(u)|\zeta(r,u)-1|\nu(\dif u)\dif r\bigg).
\de
Moreover, from (\ref{l2zetalimi}), it follows that
\ce
\lim\limits_{\d\rightarrow 0}\sup\limits_{\e}\sup\limits_{|s-t|<\d}|\Psi_\e(t)-\Psi_\e(s)|=0.
\de
Now, collecting the above estimates, by the Ascoli-Arzel\'a lemma we conclude that $\{\Psi_\e, 0<\e<1\}$ is relatively compact in $C([0, T], \mR^d)$.

Next, we observe that
\ce
\int_0^T\int_{\mU}|f(X^{\xi,\zeta}_r(x_0),u)|^2\nu(\dif u)\dif r&\leq&\int_0^T\int_{\mU}(1+|X^{\xi,\zeta}_r(x_0)|)^2L_2^2(u)\nu(\dif u)\dif r\\
&\leq&(1+\sup\limits_{r\in[0,T]}|X^{\xi,\zeta}_r(x_0)|)^2 T\int_{\mU}L_2^2(u)\nu(\dif u)\\
&\leq&C(1+|x_0|^2)T\int_{\mU}L_2^2(u)\nu(\dif u),
\de
and for any $B\subset\sB([0,T]\times\mU)$ with $(Leb_T\times\nu)(B)<\infty$,
\ce
\int_{B}\exp\{|f(X^{\xi,\zeta}_r(x_0),u)|\}\nu(\dif u)\dif r&\leq& \int_{B}\exp\{C(1+|x_0|)L_2(u)\}\nu(\dif u)\dif r\\
&\leq& \int_{B\cap\{L_2(u)>C(1+|x_0|)/\g_2\}}\exp\{C(1+|x_0|)L_2(u)\}\nu(\dif u)\dif r\\
&&+ \int_{B\cap\{L_2(u)\leq C(1+|x_0|)/\g_2\}}\exp\{C(1+|x_0|)L_2(u)\}\nu(\dif u)\dif r\\
&\leq& \int_{B\cap\{L_2(u)>C(1+|x_0|)/\g_2\}}\exp\{\g_2L^2_2(u)\}\nu(\dif u)\dif r\\
&&+\exp\{C(1+|x_0|)^2/\g_2\}\times (Leb_T\times\nu)(B)\\
&<&\infty,
\de
where $({\bf H}^{3}_{f})$ is used. So, by \cite[Lemma 3.11]{bcd}, it holds that for any $t\in[0,T]$
\ce
\lim\limits_{\e\rightarrow 0}|\Psi_\e(t)|=0, \quad \mP\mbox{-a.s.},
\de
which yields that
\ce
\lim\limits_{\e\rightarrow 0}\sup\limits_{t\in[0,T]}|\Psi_\e(t)|=0, \quad \mP\mbox{-a.s.}.
\de
The proof is complete.
\end{proof}

\bp\label{xexiezetaealsu}
Suppose that the assumptions of Theorem \ref{fwuldpth} hold. Assume that $\left(\xi^\e, \zeta^\e\right), (\xi, \zeta) \in \mathcal{A}^N$, and $\e\rightarrow 0$, $\left(\xi^\e, \zeta^\e\right)$ converges almost surely to $(\xi, \zeta)$. Then it holds that $X^{\e,\xi^\e, \zeta^\e}(x_0)$ converges in probability to $X^{\xi,\zeta}(x_0)$ as $\e\rightarrow 0$.
\ep
\begin{proof}
First of all, our aim is to prove that for $\l^\e\in {\bf \Lambda}$ with $\lim\limits_{\e\rightarrow 0}\sup\limits_{t\in[0,T]}|\l^\e(t)-t|=0$, and any $\d>0$,
\ce
\lim\limits_{\e\rightarrow 0}\mP\left\{\sup\limits_{t\in[0,T]}|X^{\e,\xi^\e, \zeta^\e}_t(x_0)-X^{\xi, \zeta}_{\l^\e(t)}(x_0)|>\d\right\}=0.
\de
Note that 
\ce
\mP\left\{\sup\limits_{t\in[0,T]}|X^{\e,\xi^\e, \zeta^\e}_t(x_0)-X^{\xi, \zeta}_{\l^\e(t)}(x_0)|>\d\right\}&\leq& \mP\left\{\sup\limits_{t\in[0,T]}|X^{\e,\xi^\e, \zeta^\e}_t(x_0)-X^{\xi, \zeta}_{t}(x_0)|>\d/2\right\}\\
&&+\mP\left\{\sup\limits_{t\in[0,T]}|X^{\xi, \zeta}_{t}(x_0)-X^{\xi, \zeta}_{\l^\e(t)}(x_0)|>\d/2\right\}.
\de
And by the similar deduction to that in the second step of the proof for Lemma \ref{xex1x2es}, it holds that
\ce
\lim\limits_{\e\rightarrow 0}\mP\left\{\sup\limits_{t\in[0,T]}|X^{\xi, \zeta}_{t}(x_0)-X^{\xi, \zeta}_{\l^\e(t)}(x_0)|>\d/2\right\}=0.
\de
Thus, we only need to show that
\be
\lim\limits_{\e\rightarrow 0}\mP\left\{\sup\limits_{t\in[0,T]}|X^{\e,\xi^\e, \zeta^\e}_t(x_0)-X^{\xi, \zeta}_{t}(x_0)|>\d/2\right\}=0.
\label{xetxt}
\ee
Our proof is divided into three steps. 

{\bf Step 1.} Set $Z^\e_t=X^{\e,\xi^\e, \zeta^\e}_t(x_0)-X^{\xi,\zeta}_t(x_0)$, and we estimate $\mE\sup\limits_{t\in[0,T]}|Z^\e_t|^2$.

First of all, by the It\^o formula, it holds that for any $t\in[0,T]$
\ce
&&|Z^\e_t|^2\\
&=&-2\int_0^t\<Z^\e_s,\dif K^{\e,\xi^\e, \zeta^\e}_s(x_0)-\dif K^{\xi,\zeta}_s(x_0)\>+2\int_0^t\<Z^\e_s, b(X^{\e,\xi^\e, \zeta^\e}_s(x_0))-b(X^{\xi,\zeta}_s(x_0))\>\dif s\\
&&+2\int_0^t\<Z^\e_s, \s(X^{\e,\xi^\e, \zeta^\e}_s(x_0))\xi^\e(s)-\s(X^{\xi,\zeta}_s(x_0))\xi(s)\>\dif s+\e\int_0^t\|\s(X^{\e,\xi^\e, \zeta^\e}_s(x_0))\|^2\dif s\\
&&+\int_0^t\int_{\mU}\<Z^\e_s, f(X^{\e,\xi^\e, \zeta^\e}_s(x_0),u)(\zeta^\e(s,u)-1)-f(X^{\xi,\zeta}_s(x_0),u)(\zeta(s,u)-1)\>\nu(\dif u)\dif s\\
&&+2\sqrt \e\int_0^t\<Z^\e_s,\s(X^{\e,\xi^\e, \zeta^\e}_s(x_0))\dif W_s\>\\
&&+\int_0^t\int_{\mU}[|Z^\e_{s-}+\e f(X^{\e,\xi^\e, \zeta^\e}_{s-}(x_0),u)|^2-|Z^\e_{s-}|^2]{\tilde N}^{\e^{-1}\zeta^\e}(\dif s\dif u)\\
&&+\e\int_0^t\int_{\mU}|f(X^{\e,\xi^\e, \zeta^\e}_{s-}(x_0),u)|^2\zeta^\e(s,u)\nu(\dif u)\dif s.
\de
Lemma \ref{equi}, $({\bf H}^{1}_{b,\s})$ and $({\bf H}^{2}_{f})$ furthermore yield that
\ce
&&\sup\limits_{s\in[0,t]}|Z^\e_s|^2\no\\
&\leq&\e(1+\sup\limits_{s\in[0,T]}|X^{\e,\xi^\e, \zeta^\e}_s(x_0)|^2)\left(CT+\int_0^T\int_{\mU}L_2^2(u)\zeta^\e(s,u)\nu(\dif u)\dif s\right)\no\\
&&+\int_0^t\(2L_1+2L_1|\xi^\e(r)|+\int_{\mU}L_2(u)|\zeta^\e(r,u)-1|\nu(\dif u)\)\sup\limits_{s\in[0,r]}|Z^\e_s|^2\dif r\no\\
&&+2\sup\limits_{s\in[0,T]}\left|\int_0^s\<Z^\e_r, \s(X^{\xi,\zeta}_r(x_0))(\xi^\e(r)-\xi(r))\>\dif r\right|\no\\
&&+\sup\limits_{s\in[0,T]}\left|\int_0^s\int_{\mU}\<Z^\e_r, f(X^{\xi,\zeta}_r(x_0),u)\((\zeta^\e(r,u)-1)-(\zeta(r,u)-1)\)\>\nu(\dif u)\dif r\right|\no\\
&&+2\sqrt \e\sup\limits_{s\in[0,T]}\left|\int_0^s\<Z^\e_r,\s(X^{\e,\xi^\e, \zeta^\e}_r(x_0))\dif W_r\>\right|\no\\
&&+\sup\limits_{s\in[0,T]}\left|\int_0^s\int_{\mU}[|Z^\e_{r-}+\e f(X^{\e,\xi^\e, \zeta^\e}_{r-}(x_0),u)|^2-|Z^\e_{r-}|^2]{\tilde N}^{\e^{-1}\zeta^\e}(\dif r\dif u)\right|.
\de
So, the Gronwall inequality implies that
\ce
&&\sup\limits_{s\in[0,t]}|Z^\e_s|^2\\
&\leq&\e C(1+\sup\limits_{s\in[0,T]}|X^{\e,\xi^\e, \zeta^\e}_s(x_0)|^2)\left(CT+M\right)\no\\
&&+2C\sup\limits_{s\in[0,T]}\left|\int_0^s\<Z^\e_r, \s(X^{\xi,\zeta}_r(x_0))(\xi^\e(r)-\xi(r))\>\dif r\right|\no\\
&&+C\sup\limits_{s\in[0,T]}\left|\int_0^s\int_{\mU}\<Z^\e_r, f(X^{\xi,\zeta}_r(x_0),u)\((\zeta^\e(r,u)-1)-(\zeta(r,u)-1)\)\>\nu(\dif u)\dif r\right|\no\\
&&+2C\sqrt \e\sup\limits_{s\in[0,T]}\left|\int_0^s\<Z^\e_r,\s(X^{\e,\xi^\e, \zeta^\e}_r(x_0))\dif W_r\>\right|\no\\
&&+C\sup\limits_{s\in[0,T]}\left|\int_0^s\int_{\mU}[|Z^\e_{r-}+\e f(X^{\e,\xi^\e, \zeta^\e}_{r-}(x_0),u)|^2-|Z^\e_{r-}|^2]{\tilde N}^{\e^{-1}\zeta^\e}(\dif r\dif u)\right|,
\de
which together with (\ref{xexezeboun}) yields that
\be
&&\mE\sup\limits_{s\in[0,T]}|Z^\e_s|^2\no\\
&\leq& C\e(1+|x_0|^2)+C\mE\sup\limits_{s\in[0,T]}\left|\int_0^s\<Z^\e_r, \s(X^{\xi,\zeta}_r(x_0))(\xi^\e(r)-\xi(r))\>\dif r\right|\no\\
&&+C\mE\sup\limits_{s\in[0,T]}\left|\int_0^s\int_{\mU}\<Z^\e_r, f(X^{\xi,\zeta}_r(x_0),u)\((\zeta^\e(r,u)-1)-(\zeta(r,u)-1)\)\>\nu(\dif u)\dif r\right|\no\\
&&+C\sqrt\e\mE\sup\limits_{s\in[0,T]}\left|\int_0^s\<Z^\e_r,\s(X^{\e,\xi^\e, \zeta^\e}_r(x_0))\dif W_r\>\right|\no\\
&&+C\mE\sup\limits_{s\in[0,T]}\left|\int_0^s\int_{\mU}[|Z^\e_{r-}+\e f(X^{\e,\xi^\e, \zeta^\e}_{r-}(x_0),u)|^2-|Z^\e_{r-}|^2]{\tilde N}^{\e^{-1}\zeta^\e}(\dif r\dif u)\right|\no\\
&=:&C\e(1+|x_0|^2)+I_1+I_2+I_3+I_4.
\label{i1i2i3i4}
\ee

We first estimate $I_3$. From the BDG inequality, it follows that
\be
I_3&\leq& C\sqrt\e\mE\left(\int_0^T|Z^\e_r|^2\|\s(X^{\e,\xi^\e, \zeta^\e}_r(x_0))\|^2\dif r\right)^{1/2}\no\\
&\leq&C\sqrt\e\mE\sup\limits_{r\in[0,T]}|Z^\e_r|\left(\int_0^T\|\s(X^{\e,\xi^\e, \zeta^\e}_r(x_0))\|^2\dif r\right)^{1/2}\no\\
&\leq&\frac{1}{4}\mE\sup\limits_{r\in[0,T]}|Z^\e_r|^2+C\e\int_0^T(1+\mE\sup\limits_{r\in[0,T]}|X^{\e,\xi^\e, \zeta^\e}_r(x_0)|^2)\dif r\no\\
&\leq&\frac{1}{4}\mE\sup\limits_{r\in[0,T]}|Z^\e_r|^2+C\e(1+|x_0|^2).
\label{i3}
\ee

For $I_4$, the BDG inequality and the mean value theorem imply that
\be
I_4&\leq& C\mE\left(\int_0^T\int_{\mU}[|Z^\e_{r-}+\e f(X^{\e,\xi^\e, \zeta^\e}_{r-}(x_0),u)|^2-|Z^\e_{r-}|^2]^2 N^{\e^{-1}\zeta^\e}(\dif r\dif u)\right)^{1/2}\no\\
&\leq&C\e\mE\left(\int_0^T\int_{\mU}|Z^\e_{r-}|^2|f(X^{\e,\xi^\e, \zeta^\e}_{r-}(x_0),u)|^2N^{\e^{-1}\zeta^\e}(\dif r\dif u)\right)^{1/2}\no\\
&&+C\e\mE\left(\int_0^T\int_{\mU}|f(X^{\e,\xi^\e, \zeta^\e}_{r-}(x_0),u)|^4N^{\e^{-1}\zeta^\e}(\dif r\dif u)\right)^{1/2}\no\\
&\leq&C\e\mE\sup\limits_{r\in[0,T]}|Z^\e_{r}|\left(\int_0^T\int_{\mU}(1+|X^{\e,\xi^\e, \zeta^\e}_{r-}(x_0)|^2)L_2^2(u)N^{\e^{-1}\zeta^\e}(\dif r\dif u)\right)^{1/2}\no\\
&&+C\e\mE\left(\int_0^T\int_{\mU}(1+|X^{\e,\xi^\e, \zeta^\e}_{r-}(x_0)|^4)L_2^2(u)N^{\e^{-1}\zeta^\e}(\dif r\dif u)\right)^{1/2}\no\\
&\leq&\frac{1}{4}\mE\sup\limits_{r\in[0,T]}|Z^\e_{r}|^2+C\e\mE\int_0^T\int_{\mU}(1+|X^{\e,\xi^\e, \zeta^\e}_{r}(x_0)|^2)L_2^2(u)\zeta^\e(r,u)\nu(\dif u)\dif r\no\\
&&+C\e\mE(1+\sup\limits_{r\in[0,T]}|X^{\e,\xi^\e, \zeta^\e}_{r}(x_0)|^2)\no\\
&&+C\e\mE\int_0^T\int_{\mU}(1+|X^{\e,\xi^\e, \zeta^\e}_{r}(x_0)|^2)L_2^2(u)\zeta^\e(r,u)\nu(\dif u)\dif r\no\\
&\leq&\frac{1}{4}\mE\sup\limits_{r\in[0,T]}|Z^\e_{r}|^2+CM\e\mE(1+\sup\limits_{r\in[0,T]}|X^{\e,\xi^\e, \zeta^\e}_{r}(x_0)|^2)+C\e(1+|x_0|^2)\no\\
&\leq&\frac{1}{4}\mE\sup\limits_{r\in[0,T]}|Z^\e_{r}|^2+C\e(1+|x_0|^2),
\label{i4}
\ee
where we use (\ref{l22zetaboun}) and (\ref{xexezeboun}). 

Finally, combining (\ref{i3}), (\ref{i4}) with (\ref{i1i2i3i4}), we conclude that
\be
\mE\sup\limits_{s\in[0,T]}|Z^\e_s|^2\leq C\e(1+|x_0|^2)+2I_1+2I_2.
\label{i1i2}
\ee
In the following step, we estimate $I_1$ and $I_2$.

{\bf Step 2.} We investigate the limits of $I_1$ and $I_2$ as $\e\rightarrow 0$.

First of all, we deal with $I_1$. Notice that
\ce
\int_0^s\<Z^\e_r, \s(X^{\xi,\zeta}_r(x_0))(\xi^\e(r)-\xi(r))\>\dif r=\int_0^s\<Z^\e_r, \dif \Sigma_\e(r)\>,
\de
where $\Sigma_\e(r)$ is defined in Lemma \ref{zeezecon}. Thus, by the It\^o formula obtain that
\be
&&\int_0^s\<Z^\e_r, \s(X^{\xi,\zeta}_r(x_0))(\xi^\e(r)-\xi(r))\>\dif r\no\\
&=&\<Z^\e_s,\Sigma_\e(s)\>+\int_0^s\<\Sigma_\e(r),\dif K^{\e,\xi^\e, \zeta^\e}_{r}(x_0)-\dif K^{\xi, \zeta}_{r}(x_0)\>\no\\
&&-\int_0^s\<\Sigma_\e(r),b(X^{\e,\xi^\e, \zeta^\e}_{r}(x_0))-b(X^{\xi, \zeta}_{r}(x_0))\>\dif r\no\\
&&-\int_0^s\<\Sigma_\e(r),\s(X^{\e,\xi^\e, \zeta^\e}_{r}(x_0))\xi^\e(r)-\s(X^{\xi, \zeta}_{r}(x_0))\xi(r)\>\dif r\no\\
&&-\sqrt \e\int_0^s\<\Sigma_\e(r),\s(X^{\e,\xi^\e, \zeta^\e}_{r}(x_0))\dif W_r\>\no\\
&&-\int_0^s\int_{\mU}\<\Sigma_\e(r),f(X^{\e,\xi^\e, \zeta^\e}_r(x_0),u)(\zeta^\e(r,u)-1)\no\\
&&\qquad\qquad -f(X^{\xi,\zeta}_r(x_0),u)(\zeta(r,u)-1)\>\nu(\dif u)\dif r\no\\
&&-\e\int_0^s\int_{\mU}\<\Sigma_\e(r), f(X^{\e,\xi^\e, \zeta^\e}_{r-}(x_0),u)\>{\tilde N}^{\e^{-1}\zeta^\e}(\dif r\dif u)\no\\
&=:&I_{11}+I_{12}+I_{13}+I_{14}+I_{15}+I_{16}+I_{17}.
\label{i11234567}
\ee
We deal with them one by one. 

For $I_{11}$, it is easy to see that
\ce
\sup\limits_{s\in[0,T]}|I_{11}(s)|\leq\sup\limits_{s\in[0,T]}|Z^\e_s|\sup\limits_{s\in[0,T]}|\Sigma_\e(s)|.
\de
On the one hand, (\ref{gelimi}) implies that
$$
\lim\limits_{\e\rightarrow 0}\sup\limits_{s\in[0,T]}|I_{11}(s)|=0.
$$
On the other hand, by the linear growth of $\s$ and (\ref{xxizeboun}), it holds that
\ce
\sup\limits_{s\in[0,T]}|\Sigma_\e(s)|&\leq& \left(\int_0^T\|\s(X^{\xi,\zeta}_r(x_0))\|(|\xi^\e(r)|+|\xi(r)|)\dif r\right)\\
&\leq&\sqrt 2\left(\int_0^T\|\s(X^{\xi,\zeta}_r(x_0)\|^2\dif r\right)^{1/2}\left(\int_0^T(|\xi^\e(r)|^2+|\xi(r)|^2)\dif r\right)^{1/2}\\
&\leq&C(1+|x_0|),
\de
which together with (\ref{xexezeboun}) and (\ref{xxizeboun}) yields that
\ce
\mE\sup\limits_{s\in[0,T]}|Z^\e_s|\sup\limits_{s\in[0,T]}|\Sigma_\e(s)|&\leq& C(1+|x_0|)\mE\sup\limits_{s\in[0,T]}|Z^\e_s|\\
&\leq& C(1+|x_0|)\mE\left(\sup\limits_{s\in[0,T]}|X^{\e,\xi^\e, \zeta^\e}_{s}(x_0)|+\sup\limits_{s\in[0,T]}|X^{\xi,\zeta}_s(x_0)|\right)\\
&\leq& C(1+|x_0|^2).
\de
Now, collecting the above deduction, by the dominated convergence theorem we obtain that
\be
\lim\limits_{\e\rightarrow 0}\mE\sup\limits_{s\in[0,T]}|I_{11}(s)|=0.
\label{i11}
\ee

Let us treat $I_{12}$. Some computations present that
\ce
\sup\limits_{s\in[0,T]}|I_{12}(s)|\leq \sup\limits_{s\in[0,T]}|\Sigma_\e(s)|(|K^{\e,\xi^\e, \zeta^\e}(x_0)|_0^T+|K^{\xi, \zeta}(x_0)|_0^T).
\de
And (\ref{gelimi}) implies that 
$$
\lim\limits_{\e\rightarrow0}\sup\limits_{s\in[0,T]}|I_{12}(s)|=0,
$$
which together with (\ref{kexezeboun}), (\ref{kxizeboun}) and the dominated convergence theorem yields that
\be
\lim\limits_{\e\rightarrow0}\mE\sup\limits_{s\in[0,T]}|I_{12}(s)|=0.
\label{i12}
\ee
By the similar deduction to that for (\ref{i12}), we can obtain that 
\be
\lim\limits_{\e\rightarrow0}\mE\sup\limits_{s\in[0,T]}(|I_{13}(s)|+|I_{14}(s)|+|I_{16}(s)|)=0.
\label{i1346}
\ee

For $I_{15}$, from the BDG inequality, it follows that
\ce
\mE\sup\limits_{s\in[0,T]}|I_{15}(s)|&\leq& C\sqrt \e\mE\left(\int_0^T|\Sigma_\e(r)|^2\|\s(X^{\e,\xi^\e, \zeta^\e}_{r}(x_0))\|^2\dif r\right)^{1/2}\\
&\leq&C\sqrt \e\mE\sup\limits_{r\in[0,T]}|\Sigma_\e(r)|^2+C\sqrt \e\int_0^T(1+\mE|X^{\e,\xi^\e, \zeta^\e}_{r}(x_0)|^2)\dif r\\
&\leq&C\sqrt \e(1+|x_0|^2),
\de
which yields that
\be
\lim\limits_{\e\rightarrow0}\mE\sup\limits_{s\in[0,T]}|I_{15}(s)|=0.
\label{i15}
\ee
Then based on the similar deduction to that for (\ref{i15}), it holds that
\be
\lim\limits_{\e\rightarrow0}\mE\sup\limits_{s\in[0,T]}|I_{17}(s)|=0.
\label{i17}
\ee

At present, collecting all the above deduction, we have that
\be
\lim\limits_{\e\rightarrow 0}I_1=0.
\label{i1}
\ee

Finally, we calculate the limit of $I_2$ as $\e\rightarrow 0$. By (\ref{psilimi}), the similar deduction to that for (\ref{i1}) yields that
\be
\lim\limits_{\e\rightarrow 0}I_2=0.
\label{i2}
\ee

{\bf Step 3.} We prove (\ref{xetxt}).

Combining (\ref{i1i2}) with (\ref{i1}) and (\ref{i2}), we obtain that
\be
\lim\limits_{\e\rightarrow 0}\mE\sup\limits_{s\in[0,T]}|Z^\e_s|^2=0,
\label{ze0}
\ee
which together with the Chebyshev inequality yields that
\ce
\lim\limits_{\e\rightarrow 0}\mP\left\{\sup\limits_{s\in[0,T]}|Z^\e_s|>\d/2\right\}=0.
\de
The proof is complete.
\end{proof}

\subsection{Proof of Theorem \ref{fwuldpth}} 

In this subsection, we verify Condition \ref{cond} and prove Theorem \ref{fwuldpth}.

{\bf Verification of Condition \ref{cond} $(i)$.} Let $\left(\xi^\e, \zeta^\e\right), (\xi, \zeta) \in \mathcal{A}^N, x_0^\e, x_0\in\overline{\cD(A)}$ be such that, as $\e\rightarrow 0, \left(\xi^\e, \zeta^\e\right)$ converges in distribution to $(\xi, \zeta)$ and $x_0^\e\rightarrow x_0$. Then since
\ce
&&\mathcal{G}^\e\left(x^\e_0, \sqrt{\e} W+\int_0^{\cdot} \xi^\e(s) \dif s, \e N^{\e^{-1} \zeta^\e}\right)=X^{\e,\xi^\e, \zeta^\e}(x^\e_0),\\
&&\mathcal{G}^0\left(x_0, \int_0^{\cdot} \xi(s) \dif s, \nu_T^{\zeta}\right)=X^{\xi, \zeta}(x_0),
\de
our goal is to prove that $X^{\e,\xi^\e, \zeta^\e}(x^\e_0)$ converges in distribution to $X^{\xi, \zeta}(x_0)$. 

Next, note that as $\e\rightarrow 0, \left(\xi^\e, \zeta^\e\right)$ converges in distribution to $(\xi, \zeta)$. By the Skorohod representation theorem, there exists a probability space $(\check{\Omega}, \check{\mathscr{F}}, \check{\mP})$, and random variables $\{(\check{\xi}^\e, \check{\zeta}^\e), 0<\e<1\}$, $(\check{\xi},\check{\zeta})$, an $l$-dimensional Brownian motion $(\check{W}_t)$, a Poisson random measure $\check{N}^{\e^{-1}\check{\zeta}^\e}(\dif t\dif u)$ with the intensity $\e^{-1}\check{\zeta}^\e(t,u)\dif t\nu(\dif u)$ such that 

$(i)$ $\sL_{(\check{\xi}^\e,\check{\zeta}^\e,\check{W},\check{N}^{\e^{-1}\check{\zeta}^\e})}=\sL_{(\xi^\e,\zeta^\e,W, N^{\e^{-1}\zeta^\e})}$ and $\sL_{(\check{\xi},\check{\zeta})}=\sL_{(\xi, \zeta)}$;

$(ii)$ $\left(\check{\xi}^\e,\check{\zeta}^\e\right)$ converges almost surely to $(\check{\xi},\check{\zeta})$ as $\e\rightarrow 0$.

Consider the following equation
\ce\left\{\begin{array}{ll}
\dif \check{X}^{\e,\check{\xi}^\e,\check{\zeta}^\e}_t\in -A(\check{X}^{\e,\check{\xi}^\e,\check{\zeta}^\e}_t)\dif t+b(\check{X}^{\e,\check{\xi}^\e,\check{\zeta}^\e}_t)\dif t+\sigma(\check{X}^{\e,\check{\xi}^\e,\check{\zeta}^\e}_t)\check{\xi}^\e(t)\dif t\\
\qquad\qquad\qquad+\sqrt \e\sigma(\check{X}^{\e,\check{\xi}^\e,\check{\zeta}^\e}_t)\dif \check{W}_t+\int_{\mU}f(\check{X}^{\e,\check{\xi}^\e,\check{\zeta}^\e}_t,u)(\check{\zeta}^\e(t,u)-1)\nu(\dif u)\dif t\\
\qquad\qquad\qquad+\e\int_{\mU}f(\check{X}^{\e,\check{\xi}^\e,\check{\zeta}^\e}_{t-},u)\tilde{\check{N}}^{\e^{-1}\check{\zeta}^\e}(\dif t\dif u),\\
\check{X}^{\e,\check{\xi}^\e,\check{\zeta}^\e}_0=x^\e_0\in \overline{\cD(A)}, \quad 0\leq t\leq T,
\end{array}
\right.
\de
where $\tilde{\check{N}}^{\e^{-1}\check{\zeta}^\e}(\dif t\dif u):=\check{N}^{\e^{-1}\check{\zeta}^\e}(\dif t\dif u)-\e^{-1}\check{\zeta}^\e(t,u)\nu(\dif u)\dif t$. 
Under $(\mathbf{H}_{A})$, $({\bf H}^{1}_{b,\s})$, $({\bf H}^{1}_{f})$ and $({\bf H}^{2}_{f})$, by Theorem 3.1 in \cite{gw}, the above equation has a unique strong solution $(\check{X}^{\e,\check{\xi}^\e,\check{\zeta}^\e}(x^\e_0), \check{K}^{\e,\check{\xi}^\e,\check{\zeta}^\e}(x^\e_0))$, and $\check{X}^{\e,\check{\xi}^\e,\check{\zeta}^\e}(x^\e_0)$ and $X^{\e, \xi^\e, \zeta^\e}(x^\e_0)$ have the same distribution, that is,
\ce
\check{X}^{\e,\check{\xi}^\e,\check{\zeta}^\e}(x^\e_0)\overset{d}{=}X^{\e, \xi^\e, \zeta^\e}(x^\e_0).
\de
By Proposition \ref{xex1x2es} and \ref{xexiezetaealsu}, we know that $\check{X}^{\e,\check{\xi}^\e,\check{\zeta}^\e}(x^\e_0)$ converges in probability to $\check{X}^{\check{\xi}, \check{\zeta}}(x_0)$ and then $\check{X}^{\e,\check{\xi}^\e,\check{\zeta}^\e}(x_0)$ converges in distribution to $\check{X}^{\check{\xi}, \check{\zeta}}(x_0)$, where $(\check{X}^{\check{\xi}, \check{\zeta}}(x_0),\check{K}^{\check{\xi}, \check{\zeta}}(x_0))$ uniquely solves the following equation
\ce\left\{\begin{array}{ll}
\dif \check{X}^{\check{\xi}, \check{\zeta}}_t\in -A(\check{X}^{\check{\xi}, \check{\zeta}}_t)\dif t+b(\check{X}^{\check{\xi}, \check{\zeta}}_t)\dif t+\sigma(\check{X}^{\check{\xi}, \check{\zeta}}_t)\check{\xi}(t)\dif t\\
\qquad\qquad\qquad\qquad+\int_{\mU}f(\check{X}^{\check{\xi}, \check{\zeta}}_t,u)(\check{\zeta}(t,u)-1)\nu(\dif u)\dif t,\\
\check{X}^{\check{\xi}, \check{\zeta}}_0=x_0\in \overline{\cD(A)}, \quad 0\leq t\leq T.
\end{array}
\right.
\de
Note that 
$$
\check{X}^{\check{\xi}, \check{\zeta}}(x_0)\overset{d}{=}X^{\xi, \zeta}(x_0).
$$
Thus, the above deduction implies that $X^{\e,\xi^\e, \zeta^\e}(x^\e_0)$ converges in distribution to $X^{\xi, \zeta}(x_0)$.

{\bf Verification of Condition \ref{cond} $(ii)$.} Proposition \ref{xex1x2es} and \ref{xexiezetaealsu} directly imply Condition \ref{cond} $(ii)$.

\vspace{3mm}

It is the position to prove Theorem \ref{fwuldpth}.

{\bf Proof of Theorem \ref{fwuldpth}.} First of all, by the above deduction we know that Condition \ref{cond} holds.

Next, we prove that for all $\phi \in D([0,T],\overline{\cD(A)})$, $\overline{\cD(A)}\ni x_0 \mapsto \Lambda_{x_0}(\phi)\in[0, \infty]$ is a lower semicontinuous mapping. So, based on the equivalence of the lower semicontinuity 
and the close property, it is sufficient to show that for any $M>0$, $\Gamma_M:=\{x_0\in\overline{\cD(A)}: \Lambda_{x_0}(\phi)\leq M\}$ is close. Assume that $\{x_0^n, n\in\mN\}\subset\Gamma_M$ and 
$x_0^n\rightarrow x_0\in\overline{\cD(A)}$ as $n\rightarrow\infty$. We prove $x_0\in\Gamma_M$.

Indeed, since $\Lambda_{x^n_0}(\phi)\leq M$, there exists a $h\in {\bf D}^{M+1/n}_1, g\in {\bf D}^{(M+1/n)/2}_2$ such that $\phi=X^{h,g}(x_0^n)$. By the similar deduction to that for Proposition \ref{xex1x2es}, it holds that $X^{h,g}(x_0^n)$ converges to $X^{h,g}(x_0)$ as $n\rightarrow\infty$. Note that $h\in {\bf D}^{M+1/n}_1, g\in {\bf D}^{(M+1/n)/2}_2$. Thus, $\Lambda_{x_0}(\phi)\leq M+1/n$. Letting $n\rightarrow\infty$, we obtain that $\Lambda_{x_0}(\phi)\leq M$ and $x_0\in\Gamma_M$. That is, $\Gamma_M$ is close.

Finally, by Theorem \ref{ulpjuth}, we conclude that Theorem \ref{fwuldpth} holds.

\subsection{Proof of Theorem \ref{dzuldpth}}

In this section, we prove Theorem \ref{dzuldpth}. We prepare a key lemma.

\bl\label{haus}
Assume that $(\mathbf{H}_{A})$, $({\bf H}^{1}_{b,\s})$, $({\bf H}^{1}_{f})$-$({\bf H}^{3}_{f})$ hold and $x_{0,n}\rightarrow x_0$ in $\overline{\cD(A)}$ as $n\rightarrow\infty$. Then for any $M\geq 0$
\ce
\lim\limits_{n\rightarrow\infty}\max\left\{\sup\limits_{\phi\in\Phi_{x_0}(M)}\rho_{D([0,T],\overline{\cD(A)})}(\phi, \Phi_{x_{0,n}}(M)), \sup\limits_{\phi\in\Phi_{x_{0,n}}(M)}\rho_{D([0,T],\overline{\cD(A)})}(\phi, \Phi_{x_0}(M))\right\}=0.
\de
\el
\begin{proof}
First of all, if $\phi\in\Phi_{x_0}(M)$, $\Lambda_{x_0}(\phi)\leq M$. By the definition of $\Lambda_{x_0}(\phi)$, there exists a $(h,g) \in\mathbf{D}$ such that $\phi=X^{h,g}(x_0)$ and 
\ce
\frac{1}{2}\int_0^T|h(t)|^2\dif t+\int_0^T\int_{\mU}\ell\(g(t,u)\)\nu(\dif u)\dif t\leq M,
\de
which also implies that $\Lambda_{x_{0,n}}(X^{h,g}(x_{0,n}))\leq M$ and $X^{h,g}(x_{0,n})\in \Phi_{x_{0,n}}(M)$. So,
\ce
\rho_{D([0,T],\overline{\cD(A)})}(\phi, \Phi_{x_{0,n}}(M))\leq \rho_{D([0,T],\overline{\cD(A)})}(\phi, X^{h,g}(x_{0,n})).
\de
Besides, note that $X^{h,g}(x_0)$ is continuous in $t$. Thus, by the similar deduction to that for (\ref{zexiexeoxo}), it holds that
\ce
\rho_{D([0,T],\overline{\cD(A)})}(\phi, \Phi_{x_{0,n}}(M))\leq\sup\limits_{t\in[0,T]}|X_t^{h,g}(x_0)-X_t^{h,g}(x_{0,n})|\leq C|x_0-x_{0,n}|,
\de
where the constant $C>0$ is independent of $n, \phi$. That is, 
\ce
\sup\limits_{\phi\in\Phi_{x_0}(M)}\rho_{D([0,T],\overline{\cD(A)})}(\phi, \Phi_{x_{0,n}}(M))\leq C|x_0-x_{0,n}|.
\de

Similarly, it holds that
\ce
\sup\limits_{\phi\in\Phi_{x_{0,n}}(M)}\rho_{D([0,T],\overline{\cD(A)})}(\phi, \Phi_{x_0}(M))\leq C|x_0-x_{0,n}|.
\de
Finally, the above deduction implies the required result.
\end{proof}

Now, we prove Theorem \ref{dzuldpth}.

{\bf Proof of Theorem \ref{dzuldpth}.} Lemma \ref{haus} and \cite[Theorem 2.7]{ms} assure that Theorem \ref{dzuldpth} holds.

\section{Proof of Theorem \ref{ldpinvameas}}\label{ldpinvameasproo}

In this section, we prove Theorem \ref{ldpinvameas}. By \cite[Section 4.4]{fw}, it is sufficient to prove three following points:

$(I)$ compactness: for any $s\geq 0$, ${\bf K}(s)=\{y\in\overline{\cD(A)}: V(y)\leq s\}$ is compact.

$(II)$ lower bounds: for any $\d>0, \t>0$ and $y^*\in\overline{\cD(A)}$, there exists a $\e_0>0$ such that
\ce
\mu^\e(y\in\overline{\cD(A)}: |y-y^*|<\d)\geq \exp\left\{-\frac{V(y^*)+\t}{\e}\right\},
\de
for all $0<\e<\e_0$.

$(III)$ upper bounds: for any $s\geq 0, \d>0, \t>0$, there exists a $\e_0>0$ such that 
\ce
\mu^\e(y\in\overline{\cD(A)}: \rho_{\overline{\cD(A)}}(y, {\bf K}(s))\geq \d)\leq \exp\left\{-\frac{s-\t}{\e}\right\},
\de
for all $0<\e<\e_0$, where 
$$
\rho_{\overline{\cD(A)}}(y, {\bf K}(s)):=\inf\limits_{z\in{\bf K}(s)}|y-z|.
$$
We prove them in three following subsections.

\subsection{Compactness.} In this subsection, we prove the compactness. And we begin with some key lemmas.

Consider Eq.(\ref{x0eq}), i.e.
\be\left\{\begin{array}{ll}
\dif X^0_t\in -A(X^0_t)\dif t+b(X^0_t)\dif t, \\
X^0_0=x_0\in \overline{\cD(A)}, \quad t\geq 0.
\end{array}
\right.
\label{x0eq1}
\ee

\bl\label{xosoexes}
Under $(\mathbf{H}_{A})$, $({\bf H}^{1}_{b,\s})$ and $({\bf H}_{b,\s,f})$, Eq.(\ref{x0eq1}) has a unique solution $(X^0(x_0),K^0(x_0))$. Moreover, it holds that for any $t\geq 0$
\be
|X^0_t(x_0)|^2\leq |x_0|^2e^{-L_3 t}.
\label{x0es}
\ee
\el
\begin{proof}
First of all, under $(\mathbf{H}_{A})$ and $({\bf H}^{1}_{b,\s})$, Eq.(\ref{x0eq1}) has a unique solution $(X^0(x_0),K^0(x_0))$. Then for any $\eta>0$, we take the Yosida approximation $A^\eta$ of the operator $A$ (cf. \cite{qg}) and consider the following equation
\ce\left\{\begin{array}{ll}
\dif X^{0,\eta}_t=-A^\eta(X^{0,\eta}_t)\dif t+b(X^{0,\eta}_t)\dif t, \\
X^{0,\eta}_0=x_0\in \overline{\cD(A)}, \quad t\geq 0.
\end{array}
\right.
\de
Since $A^\eta$ is a single-valued and Lipschitz continuous function, the above equation has a unique solution $X^{0,\eta}(x_0)$. Moreover, some computations imply that 
\ce
\dif |X^{0,\eta}_t(x_0)|^2&=&-2\<X^{0,\eta}_t(x_0), A^\eta(X^{0,\eta}_t(x_0))\>\dif t+2\<X^{0,\eta}_t(x_0), b(X^{0,\eta}_t(x_0))\>\dif t\\
&\leq&-L_3|X^{0,\eta}_t(x_0)|^2\dif t,
\de
where we use the maximal monotone property of $A^\eta$ and $({\bf H}_{b,\s,f})$. By the comparison theorem for ordinary differential equations, it holds that
\ce
|X^{0,\eta}_t(x_0)|^2\leq |x_0|^2e^{-L_3 t}.
\de
Note that $X^{0,\eta}(x_0)$ converges to $X^0(x_0)$ and the right side of the above inequality is independent of $\eta$. Thus, we obtain (\ref{x0es}).
\end{proof}

Next, we have the following result.

\bl\label{lam0t}
Assume that $(\mathbf{H}_{A})$, $({\bf H}^{1}_{b,\s})$, $({\bf H}^2_{f})$, $({\bf H}^3_{f})$ and $({\bf H}_{b,\s,f})$ hold. For any $\kappa\geq 0, r\geq 0$ and $T>0$, set
$$
S(\kappa, r, T):=\left\{\phi\in D([0,T], \overline{\cD(A)}): |\phi(0)|\leq \kappa, |\phi(T)|\geq r\right\}.
$$
Then for any $\phi\in S(\kappa, r, T)$, $\Lambda_0^T(\phi)\rightarrow\infty$ as $T\rightarrow\infty$.
\el
\begin{proof}
First of all, we establish that there exists a $T_0>0$ such that
\ce
\t:=\inf\{\Lambda_0^{T_0}(\phi): \phi\in S(\kappa, r, T_0)\}>0.
\de
If it is not true, for any $T>0$ there exists a sequence $\{\phi_n\}\subset S(\kappa, r, T)$ such that 
$$
\lim\limits_{n\rightarrow\infty}\Lambda_0^{T}(\phi_n)=0.
$$
Based on the definition of $\Lambda_0^{T}(\phi_n)$, there also exists a sequence $\{(h_n, g_n)\}\subset\mathbf{D}$ such that $\phi_n=X^{h_n,g_n}$, and 
\ce
\lim\limits_{n\rightarrow\infty}\int_0^T|h_n(t)|^2\dif t=0, \quad \lim\limits_{n\rightarrow\infty}\int_0^T\int_{\mU}\ell(g_n(t,u))\nu(\dif u)\dif t=0.
\de
On one hand, $(\phi_n, k_n)$ satisfies the following equation:
\ce
\phi_n(t)&=&\phi_n(0)- k_n(t)+\int_0^t b(\phi_n(s))\dif s+\int_0^t \sigma(\phi_n(s))h_n(s)\dif s\\
&&+\int_0^t\int_{\mU}f(\phi_n(s),u)(g_n(s,u)-1)\nu(\dif u)\dif s.
\de
On the other hand, there exists a $N\in\mN$ such that $\{(h_n, g_n)\}\subset\mathbf{D}^N$ and $(h_n,g_n)\rightarrow(0,1)$ as $n$ tends to $\infty$. So, the similar deduction to that for (\ref{zexiexeoxo}) and (\ref{ze0}) yields that $\phi_n(t)$ converges to $X^0_t(\a)$ for any $t\in[0,T]$, where $\a:=\lim\limits_{n\rightarrow\infty}\phi_n(0)$ and $|\a|\leq\kappa$. Note that $\phi_n\in S(\kappa, r, T)$ and $|\phi_n(T)|\geq r$. Thus, $|X^0_T(\a)|\geq r$. However, for this $r$, by (\ref{x0es}), there exists a $\hat{T}>0$ such that for any $T\geq \hat{T}$, $|X^0_T(\a)|<r$, which is a contradiction. Hence, $\t>0$.

Next, for any $T\geq nT_0$, we know that $\Lambda_0^{T}(\phi)\geq n\t$ and the required result holds.
\end{proof}

Now, for any $(h, g) \in\mathbf{D}$, we consider the following multivalued differential equation:
\ce\left\{\begin{array}{ll}
\dif \phi(t)\in -A(\phi(t))\dif t+b(\phi(t))\dif t+\sigma(\phi(t))h(t)\dif t\\
\qquad\qquad\qquad\qquad+\int_{\mU}f(\phi(t),u)(g(t,u)-1)\nu(\dif u)\dif t,\\
\phi(0)=0\in \overline{\cD(A)}, \quad 0\leq t\leq T.
\end{array}
\right.
\de
Under $(\mathbf{H}_{A})$, $({\bf H}^{1}_{b,\s})$ and $({\bf H}^2_{f})$, the above equation has a unique solution $(\phi, k)$. Moreover, we have the following result about $\phi$.

\bl\label{phiini0li}
Suppose that $(\mathbf{H}_{A})$, $({\bf H}^{1}_{b,\s})$, $({\bf H}^2_{f})$ and $({\bf H}^3_{f})$ hold. Then $T$ converges to $\infty$ as $|\phi(T)|$ tends to $\infty$.
\el
\begin{proof}
Assume that $T$ converges to a finite number as $|\phi(T)|$ tends to $\infty$. Then we investigate $\phi(t)$ for any $t\in[0,T]$. 

Note that $(\phi, k)$ satisfies the following equation: for $t\in[0,T]$
\ce
\phi(t)=-k(t)+\int_0^t b(\phi(s))\dif s+\int_0^t \sigma(\phi(s))h(s)\dif s+\int_0^t\int_{\mU}f(\phi(s),u)(g(s,u)-1)\nu(\dif u)\dif s.
\de
Thus, the linear growth of $b,\s, f$ and some calculations imply that 
\ce
|\phi(t)|^2&=&-2\int_0^t\<\phi(s),\dif k(s)\>+2\int_0^t\<\phi(s),b(\phi(s))\>\dif s+2\int_0^t\<\phi(s),\sigma(\phi(s))h(s)\>\dif s\\
&&+2\int_0^t\int_{\mU}\<\phi(s),f(\phi(s),u)(g(s,u)-1)\>\nu(\dif u)\dif s\\
&\leq&2C\int_0^t|\phi(s)|(1+|\phi(s)|)\dif s+2C\int_0^t|\phi(s)|(1+|\phi(s)|)|h(s)|\dif s\\
&&+2\int_0^t\int_{\mU}|\phi(s)|(1+|\phi(s)|)L_2(u)|g(s,u)-1|\nu(\dif u)\dif s,
\de
and 
\ce
1+|\phi(t)|^2\leq 1+\int_0^t\left(4C+4C|h(s)|+4\int_{\mU}L_2(u)|g(s,u)-1|\nu(\dif u)\right)(1+|\phi(s)|^2)\dif s.
\de
Since $(h, g) \in\mathbf{D}$, there exists a $N\in\mN$ such that $(h, g) \in\mathbf{D}^N$, which together with (\ref{xiboun}), (\ref{l2zetaboun}) and the Gronwall inequality yields that 
$$
|\phi(T)|\leq C<\infty.
$$
When $T$ converges to a finite number, $|\phi(T)|$ converges to a finite number, which is contradict with that $|\phi(T)|$ tends to $\infty$.

Finally, by the above deduction, we conclude that $T\rightarrow\infty$ as $|\phi(T)|\rightarrow \infty$.
\end{proof}

Having Lemma \ref{lam0t} and \ref{phiini0li} in hand, we will obtain the following result.

\bl\label{vyinf}
Assume that $(\mathbf{H}_{A})$, $({\bf H}^{1}_{b,\s})$, $({\bf H}^2_{f})$, $({\bf H}^3_{f})$ and $({\bf H}_{b,\s,f})$ hold. Then $V(y)$ converges to $\infty$ as $|y|$ tends to $\infty$.
\el
\begin{proof}
By the definition of $V(y)$, for any $\eta>0$, there exists a $\phi$ such that $\phi(0)=0, \phi(T)=y$ and $\Lambda_0^{T}(\phi)\leq V(y)+\eta$. As $|y|$ tends to $\infty$, $|\phi(T)|$ converges to $\infty$. Then by Lemma \ref{phiini0li}, 
it holds that $T\rightarrow\infty$, which together with Lemma \ref{lam0t} yields that $\Lambda_0^{T}(\phi)\rightarrow\infty$ and $V(y)$ converges to $\infty$.
\end{proof}

Now, we prove the compactness.

{\bf Proof of $(I)$.} First of all, we show that ${\bf K}(s)$ is close. Assume that $\{y_j\}\subset{\bf K}(s)$ and $y_j$ converges to $y$ in $\overline{\cD(A)}$. So, for any $\d_1>0$, there exists a $j_0\in\mN$ such that 
\ce
|y_{j_0}-y|<\frac{\d_1}{2}.
\de
Besides, note that 
\ce
V(y)=\lim\limits_{\d\rightarrow0}\inf\{\Lambda_0^T(\phi): \phi(0)=0, \phi(T)\in B(y, \d)\cap\overline{\cD(A)}, T>0\},
\de
where $B(y,\d):=\{z\in\mR^d: |z-y|\leq \d\}$ (cf. \cite{dm1}). Then, since $V(y_{j_0})\leq s$, for any $\eta>0$, there exist $\phi, T$ such that $\phi(0)=0, \phi(T)\in B(y_{j_0}, \frac{\d_1}{2})\cap\overline{\cD(A)}$ and
\ce
\Lambda_0^T(\phi)\leq s+\eta.
\de
Based on the definition of $\Lambda_0^T(\phi)$, there exist a $(h,g)\in\mathbf{D}$ such that $\phi=X^{h, g}$ and
$$
\frac{1}{2}\int_0^{T}|h(t)|^2\dif t+\int_0^{T}\int_{\mU}\ell\(g(t,u)\)\nu(\dif u)\dif t\leq s+\eta.
$$
That is, $\phi(0)=0, |\phi(T)-y|\leq \d_1$, $\Lambda_0^T(\phi)\leq s+\eta$ and $V(y)\leq s+\eta$. The arbitrary of $\eta$ assures that $y\in{\bf K}(s)$.

Next, we prove that ${\bf K}(s)$ is bounded. If ${\bf K}(s)$ is unbounded, there exists a sequence $\{y_n\}\subset\overline{\cD(A)}$ satisfying that $|y_n|\rightarrow\infty$ as $n\rightarrow\infty$, and $V(y_n)\leq s$. However, by Lemma \ref{vyinf} one can obtain that as $|y_n|\rightarrow\infty$, $V(y_n)$ converges to $\infty$. For this $s$, there exists a $n_0\in\mN$ such that for any $n>n_0$, $V(y_n)>s$. This is a contradiction. Therefore, ${\bf K}(s)$ is bounded. 

\subsection{Lower bounds.} In this subsection, we prove lower bounds. We begin with the exponential moment estimate of $X^\e(x_0)$.

\bl\label{expoesti}
Suppose that $(\mathbf{H}_{A})$, $({\bf H}^{1}_{b,\s})$, $({\bf H}^{2}_{\s})$, $({\bf H}^{1}_{f})$, $({\bf H}^{2'}_{f})$, $({\bf H}^{3}_{f})$ and $({\bf H}_{b,\s,f})$ hold. Then it holds that
\ce
\mE e^{\frac{\b}{\e}|X^\e_t(x_0)|^2}\leq e^{-\frac{L_3}{4} t} e^{\frac{\b}{\e}|x_0|^2}+2, \quad t\geq 0,
\de
where the constant $\b>0$ only depends on $L_3, L_\s, \int_{\mU}L_2^2(u)\nu(\dif u)$.
\el
\begin{proof}
First of all, note that $(X^\e(x_0), K^\e(x_0))$ satisfies the following equation
\ce
X^\e_t(x_0)&=&x_0-K^\e_t(x_0)+\int_0^t b(X^\e_s(x_0))\dif s+\sqrt \e\int_0^t\sigma(X^\e_s(x_0))\dif W_s\\
&&+\e\int_0^t\int_{\mU}f(X^\e_{s-}(x_0),u) \tilde N^{\e^{-1}}(\dif s\dif u), \quad t\geq 0.
\de
Thus, by the It\^o formula, it holds that 
\ce
|X^\e_t(x_0)|^2&=&|x_0|^2-\int_0^t 2\<X^\e_s(x_0),\dif K^\e_s(x_0)\>+\int_0^t 2\<X^\e_s(x_0),b(X^\e_s(x_0))\>\dif s\\
&&+\e\int_0^t\|\sigma(X^\e_s(x_0))\|^2\dif s+\e\int_0^t\int_{\mU}|f(X^\e_{s-}(x_0),u)|^2\nu(\dif u)\dif s\\
&&+2\sqrt \e\int_0^t\<X^\e_s(x_0),\sigma(X^\e_s(x_0))\dif W_s\>\\
&&+\int_0^t\int_{\mU}\left[|X^\e_s(x_0)+\e f(X^\e_{s-}(x_0),u)|^2-|X^\e_s(x_0)|^2\right]\tilde N^{\e^{-1}}(\dif s\dif u).
\de
For any $\d>0, \b>0$ and $0<\e<1$, applying the It\^o formula to $e^{\d t+\frac{\b}{\e}|X^\e_t(x_0)|^2}$, by Lemma \ref{equi}, $({\bf H}^{2}_{\s})$, $({\bf H}^{2'}_{f})$ and $({\bf H}_{b,\s,f})$ we obtain that
\ce
e^{\d t+\frac{\b}{\e}|X^\e_t(x_0)|^2}&\leq&e^{\frac{\b}{\e}|x_0|^2}+\d\int_0^te^{\d s+\frac{\b}{\e}|X^\e_s(x_0)|^2}\dif s-L_3\frac{\b}{\e}\int_0^te^{\d s+\frac{\b}{\e}|X^\e_s(x_0)|^2}|X^\e_s(x_0)|^2\dif s\\
&&+2\frac{\b}{\e}\sqrt \e\int_0^t e^{\d s+\frac{\b}{\e}|X^\e_s(x_0)|^2}\<X^\e_s(x_0),\sigma(X^\e_s(x_0))\dif W_s\>\\
&&+\frac{\b}{\e}\int_0^t\int_{\mU}e^{\d s+\frac{\b}{\e}|X^\e_s(x_0)|^2}\left[|X^\e_s(x_0)+\e f(X^\e_{s-}(x_0),u)|^2-|X^\e_s(x_0)|^2\right]\tilde N^{\e^{-1}}(\dif s\dif u)\\
&&+2\frac{\b^2}{\e}\left(L_\s^2+\int_{\mU}L_2^2(u)\nu(\dif u)\right)\int_0^t e^{\d s+\frac{\b}{\e}|X^\e_s(x_0)|^2}|X^\e_s(x_0)|^2\dif s\\
&&+2\b^2\left(\int_{\mU}L_2^2(u)\nu(\dif u)\right)\int_0^te^{\d s+\frac{\b}{\e}|X^\e_s(x_0)|^2}\dif s.
\de
Let $\d=\frac{L_3}{4}$, and choose $\b$ small enough such that
\ce
2\b L_\s^2\leq \frac{L_3}{4}, \quad 2\b\int_{\mU}L_2^2(u)\nu(\dif u)\leq \frac{L_3}{4}, \quad\b\leq 1.
\de
So, it holds that
\ce
e^{\frac{L_3}{4} t+\frac{\b}{\e}|X^\e_t(x_0)|^2}&\leq&e^{\frac{\b}{\e}|x_0|^2}+\frac{L_3}{2}\int_0^te^{\frac{L_3}{4} s+\frac{\b}{\e}|X^\e_s(x_0)|^2}(1-\frac{\b}{\e}|X^\e_s(x_0)|^2)\dif s\\
&&+2\frac{\b}{\e}\sqrt \e\int_0^t e^{\frac{L_3}{4} s+\frac{\b}{\e}|X^\e_s(x_0)|^2}\<X^\e_s(x_0),\sigma(X^\e_s(x_0))\dif W_s\>\\
&&+\frac{\b}{\e}\int_0^t\int_{\mU}e^{\frac{L_3}{4} s+\frac{\b}{\e}|X^\e_s(x_0)|^2}\left[|X^\e_s(x_0)+\e f(X^\e_{s-}(x_0),u)|^2-|X^\e_s(x_0)|^2\right]\tilde N^{\e^{-1}}(\dif s\dif u).
\de
Taking the expectation on two sides of the above inequality, by the fact that $e^v(1-v)\leq 1$ for any $v\geq 0$ we have that
\ce
\mE e^{\frac{L_3}{4} t+\frac{\b}{\e}|X^\e_t(x_0)|^2}\leq e^{\frac{\b}{\e}|x_0|^2}+\frac{L_3}{2}\int_0^te^{\frac{L_3}{4} s}\dif s\leq e^{\frac{\b}{\e}|x_0|^2}+2 e^{\frac{L_3}{4} t},
\de
and furthermore
\ce
\mE e^{\frac{\b}{\e}|X^\e_t(x_0)|^2}\leq e^{-\frac{L_3}{4} t} e^{\frac{\b}{\e}|x_0|^2}+2.
\de
\end{proof}

\bl\label{muees}
Suppose that $(\mathbf{H}_{A})$, $({\bf H}^{1}_{b,\s})$, $({\bf H}^{2}_{\s})$, $({\bf H}^{1}_{f})$, $({\bf H}^{2'}_{f})$, $({\bf H}^{3}_{f})$ and $({\bf H}_{b,\s,f})$ hold. Then for any $\eta>0$, there exist $r_0>0$ and $\e_0>0$ such that for all $r>r_0, 0<\e<\e_0$
\ce
\mu^\e(B^c(r)\cap\overline{\cD(A)})\leq \exp\(-\frac{\eta}{\e}\).
\de
\el
\begin{proof}
Note that $\mu^\e$ is an invariant probability measure for the Markov process $X^\e(x_0)$. Thus, by the Chebyshev inequality and Lemma \ref{expoesti}, it holds that for any $t>0$
\ce
\mu^\e(B^c(r)\cap\overline{\cD(A)})&=&\int_{\overline{\cD(A)}}\mP(X^\e_t(x_0)\in B^c(r)\cap\overline{\cD(A)})\mu^\e(\dif x_0)\\
&=&\int_{\overline{\cD(A)}}\mP(X^\e_t(x_0)\in\overline{\cD(A)}, e^{\frac{\b}{\e}|X^\e_t(x_0)|^2}>e^{\frac{\b}{\e}r^2})\mu^\e(\dif x_0)\\
&\leq&\int_{\overline{\cD(A)}}e^{-\frac{\b}{\e}r^2}\mE e^{\frac{\b}{\e}|X^\e_t(x_0)|^2}\mu^\e(\dif x_0)\\
&\leq&\int_{\overline{\cD(A)}}e^{-\frac{\b}{\e}r^2}e^{-\frac{L_3}{4} t} e^{\frac{\b}{\e}|x_0|^2}\mu^\e(\dif x_0)+2e^{-\frac{\b}{\e}r^2}.
\de
As $t\rightarrow \infty$, the Fatou lemma implies that
\ce
\mu^\e(B^c(r)\cap\overline{\cD(A)})\leq 2e^{-\frac{\b}{\e}r^2}.
\de
For any $\eta>0$, there exist $\e_0>0$ and $r_0=\sqrt{\frac{\eta}{\b}+\frac{\e_0 ln 2}{\b}}$ such that for all $r>r_0, 0<\e<\e_0$
\ce
\mu^\e(B^c(r)\cap\overline{\cD(A)})\leq \exp\(-\frac{\eta}{\e}\).
\de
The proof is complete.
\end{proof}

At present we prove the lower bounds.

{\bf Proof of $(II)$.} If there exists a $y^*\in\overline{\cD(A)}$ such that $V(y^*)=\infty$, the result obviously holds. Therefore, we assume that for any $y^*\in\overline{\cD(A)}$, $V(y^*)<\infty$. By the definition of $V(y^*)$, for any $\t>0$, there exists a $T_0>0$ and $\bar\phi\in D([0,T_0],\overline{\cD(A)})$ such that $\bar\phi(0)=0, \bar\phi(T_0)=y^*$ and 
$$
\Lambda_0^{T_0}(\bar\phi)\leq V(y^*)+\frac{\t}{2}.
$$
Based on the definition of $\Lambda_0^{T_0}(\bar\phi)$, there exists $(\bar h,\bar g)\in{\bf D}$ such that $\bar\phi=X^{\bar h,\bar g}(0)$. That is, $(\bar\phi,\bar k)$ satisfies the following equation:
\ce\left\{\begin{array}{ll}
\dif \bar\phi(t)\in -A(\bar\phi(t))\dif t+b(\bar\phi(t))\dif t+\sigma(\bar\phi(t))\bar h(t)\dif t\\
\qquad\qquad\qquad\qquad+\int_{\mU}f(\bar\phi(t),u)(\bar g(t,u)-1)\nu(\dif u)\dif t,\\
\bar\phi(0)=0\in \overline{\cD(A)}, \quad 0\leq t\leq T_0.
\end{array}
\right.
\de

Next, for any $T>0$, set
\ce
&&\hat h(t):=\left\{\begin{array}{ll} 0, \qquad\qquad t\in[0,T],\\
\bar h(t-T), \quad t\in [T, T+T_0],
\end{array}
\right.\\
&&\hat g(t,u):=\left\{\begin{array}{ll} 1, \quad\qquad\qquad (t,u)\in[0,T]\times\mU,\\
\bar g(t-T,u), \quad (t,u)\in [T, T+T_0]\times\mU,
\end{array}
\right.
\de
and we construct the following multivalued differential equation:
\ce\left\{\begin{array}{ll}
\dif \hat\phi(t)\in -A(\hat\phi(t))\dif t+b(\hat\phi(t))\dif t+\sigma(\hat\phi(t))\hat h(t)\dif t\\
\qquad\qquad\qquad\qquad+\int_{\mU}f(\hat\phi(t),u)(\hat g(t,u)-1)\nu(\dif u)\dif t,\\
\hat\phi(0)=x_0\in \overline{\cD(A)}, \quad 0\leq t\leq T+T_0.
\end{array}
\right.
\de
The above equation has a unique solution $(\hat\phi, \hat k)$. Moreover, $\hat\phi(t)=X^0_t(x_0)$ for any $t\in[0,T]$ and $\hat\phi(t)=X^{\bar h,\bar g}_{t-T}(X^0_T(x_0))$ for any $t\in[T,T+T_0]$. Put $\check\phi(t):=\hat\phi(t+T)$ for any $t\in[0,T_0]$, and $(\check\phi,\check k)$ satisfies the following equation
\ce\left\{\begin{array}{ll}
\dif \check\phi(t)\in -A(\check\phi(t))\dif t+b(\check\phi(t))\dif t+\sigma(\check\phi(t))\bar h(t)\dif t\\
\qquad\qquad\qquad\qquad+\int_{\mU}f(\check\phi(t),u)(\bar g(t,u)-1)\nu(\dif u)\dif t,\\
\check\phi(0)=X^0_T(x_0)\in \overline{\cD(A)}, \quad 0\leq t\leq T_0.
\end{array}
\right.
\de
We compute $|\check\phi(t)-\bar\phi(t)|^2$ for any $t\in[0,T_0]$. By Lemma \ref{equi}, $({\bf H}^{1}_{b,\s})$ and $({\bf H}^{2'}_{f})$, it holds that
\ce
&&|\check\phi(t)-\bar\phi(t)|^2\\
&=&|X^0_T(x_0)|^2-2\int_0^t\<\check\phi(s)-\bar\phi(s),\dif \check k(s)-\bar k(s)\>\\
&&+2\int_0^t\<\check\phi(s)-\bar\phi(s),b(\check\phi(s))-b(\bar\phi(s))\>\dif s\\
&&+2\int_0^t\<\check\phi(s)-\bar\phi(s),(\s(\check\phi(s))-\s(\bar\phi(s)))\bar h(s)\>\dif s\\
&&+2\int_0^t\int_{\mU}\<\check\phi(s)-\bar\phi(s),(f(\check\phi(s),u)-f(\bar\phi(s),u))(\bar g(s,u)-1)\>\nu(\dif u)\dif s\\
&\leq&|X^0_T(x_0)|^2+\int_0^t\left(2L_1+2L_1|\bar h(s)|+\int_{\mU}2L_2(u)|\bar g(s,u)-1|\nu(\dif u)\right)|\check\phi(s)-\bar\phi(s)|^2\dif s.
\de
Since $(\bar h,\bar g)\in{\bf D}$, there exists a $N\in\mN$ such that $(\bar h,\bar g)\in{\bf D}^N$, which together with (\ref{xiboun}), (\ref{l2zetaboun}) and the Gronwall inequality implies that
\be
|\check\phi(T_0)-\bar\phi(T_0)|^2\leq |X^0_T(x_0)|^2\exp\left\{\int_0^{T_0}\left(2L_1+2L_1|\bar h(s)|+\int_{\mU}2L_2(u)|\bar g(s,u)-1|\nu(\dif u)\right)\dif s\right\}.
\label{checkbar}
\ee
Based on this and Lemma \ref{xosoexes}, for any $\d>0$, there exists a $T_1>0$ such that for $T\geq T_1$
$$
|\check\phi(T_0)-\bar\phi(T_0)|<\d/2,
$$
which together with $\check\phi(T_0)=\hat\phi(T_0+T_1), \bar\phi(T_0)=y^*$ yields that
\ce
|X^\e_{T_1+T_0}(x_0)-y^*|&\leq&|X^\e_{T_1+T_0}(x_0)-\hat\phi(T_0+T_1)|+|\hat\phi(T_0+T_1)-y^*|\\
&=&|X^\e_{T_1+T_0}(x_0)-\hat\phi(T_0+T_1)|+|\check\phi(T_0)-\bar\phi(T_0)|\\
&<&|X^\e_{T_1+T_0}(x_0)-\hat\phi(T_0+T_1)|+\d/2.
\de

Besides, by Lemma \ref{muees}, there exists $r_0>0$ and $\e_0>0$ such that for any $0<\e<\e_0$
$$
\mu^\e(B(r_0)\cap\overline{\cD(A)})\geq \frac{1}{2}.
$$

Finally, collecting the above deduction, by Remark \ref{fwthagst} $(ii)$ we conclude that
\ce
\mu^\e(y\in\overline{\cD(A)}: |y-y^*|<\d)&=&\int_{\overline{\cD(A)}}\mP(|X^\e_{T_1+T_0}(x_0)-y^*|<\d)\mu^\e(\dif x_0)\\
&\geq&\int_{\overline{\cD(A)}}\mP(|X^\e_{T_1+T_0}(x_0)-\hat\phi(T_0+T_1)|<\d/2)\mu^\e(\dif x_0)\\
&\geq&\int_{\overline{\cD(A)}}\mP(\rho_{D([0,T_1+T_0],\overline{\cD(A)})}(X^\e(x_0),\hat\phi)<\d/2)\mu^\e(\dif x_0)\\
&\geq&\int_{B(r_0)\cap\overline{\cD(A)}}\exp\left\{-\frac{\Lambda_{x_0}(\hat\phi)+\t/2}{\e}\right\}\mu^\e(\dif x_0)\\
&=&\int_{B(r_0)\cap\overline{\cD(A)}}\exp\left\{-\frac{\Lambda_0^{T_0}(\bar\phi)+\t/2}{\e}\right\}\mu^\e(\dif x_0)\\
&\geq&\exp\left\{-\frac{V(y^*)+\t}{\e}\right\}\mu^\e(B(r_0)\cap\overline{\cD(A)})\\
&\geq&\frac{1}{2}\exp\left\{-\frac{V(y^*)+\t}{\e}\right\}.
\de
The proof is complete.

\subsection{Upper bounds.} In this subsection, we prove upper bounds. 

\bl\label{kbrks}
Suppose that $(\mathbf{H}_{A})$, $({\bf H}^{1}_{b,\s})$, $({\bf H}^{2}_{f})$, $({\bf H}^{3}_{f})$ and $({\bf H}_{b,\s,f})$ hold. Then for any $\d>0, s>0$, there exists a $R^{'}>0$ and $T'>0$ such that for any $R\leq R^{'}$ and $T\geq T'$
\ce
\{\phi(T): \phi\in {\bf K}_{B(R), 0, T}(s)\}\subset \{y\in\overline{\cD(A)}: \rho_{\overline{\cD(A)}}(y, {\bf K}(s))<\d/2\},
\de
where ${\bf K}_{B(R), 0, T}(s):=\{\phi\in D([0,T],\overline{\cD(A)}): \phi(0)\in B(R)\cap\overline{\cD(A)}, \Lambda_0^T(\phi)\leq s\}$.
\el
\begin{proof}
If this is not true, there exist $R_n\rightarrow 0, T_n\rightarrow \infty$ as $n\rightarrow\infty$ and $\phi_n\in {\bf K}_{B(R_n), 0, T_n}(s)$ such that $\rho_{\overline{\cD(A)}}(\phi_n(T_n), {\bf K}(s))\geq\d/2$.

On one hand, $\phi_n\in {\bf K}_{B(R_n), 0, T_n}(s)$ implies that $\phi_n(0)\in B(R_n)\cap\overline{\cD(A)}, \Lambda_0^{T_n}(\phi)\leq s$. Thus, there exists a $(h_n,g_n)\in \bf D$ such that $\phi_n=X^{h_n,g_n}$ and 
\be
\frac{1}{2}\int_0^{T_n}|h_n(t)|^2\dif t+\int_0^{T_n}\int_{\mU}\ell\(g_n(t,u)\)\nu(\dif u)\dif t\leq s.
\label{hngnboun}
\ee
Besides, we consider the following multivalued differential equation:
\ce\left\{\begin{array}{ll}
\dif \tilde\phi_n(t)\in -A(\tilde\phi_n(t))\dif t+b(\tilde\phi_n(t))\dif t+\sigma(\tilde\phi_n(t))h_n(t)\dif t\\
\qquad\qquad\qquad\qquad+\int_{\mU}f(\tilde\phi_n(t),u)(g_n(t,u)-1)\nu(\dif u)\dif t,\\
\tilde\phi_n(0)=0\in \overline{\cD(A)}, \quad 0\leq t\leq T_n.
\end{array}
\right.
\de
Under $(\mathbf{H}_{A})$, $({\bf H}^{1}_{b,\s})$ and $({\bf H}^2_{f})$, the above equation has a unique solution $(\tilde\phi_n, \tilde k_n)$. Moreover, by the similar deduction to that for (\ref{checkbar}), it holds that
\ce
|\tilde\phi_n(T_n)-\phi_n(T_n)|^2\leq |\phi_n(0)|^2\exp\left\{\int_0^{T_n}\left(2L_1+2L_1|h_n(t)|+\int_{\mU}2L_2(u)| g_n(t,u)-1|\nu(\dif u)\right)\dif t\right\}.
\de
Note that $|\phi_n(0)|\leq R_n$ and $R_n\rightarrow 0$. Thus there exists a $n\in\mN$ such that 
$$
|\tilde\phi_n(T_n)-\phi_n(T_n)|\leq \d/4,
$$
which together with $\rho_{\overline{\cD(A)}}(\phi_n(T_n), {\bf K}(s))\geq\d/2$ yields that
\be
\rho_{\overline{\cD(A)}}(\tilde\phi_n(T_n), {\bf K}(s))\geq\d/4.
\label{tildephi}
\ee

On the other hand, the definition of the rate function $V$ and (\ref{hngnboun}) assures that
\ce
V(\tilde\phi_n(T_n))\leq \Lambda_0^{T_n}(\tilde\phi_n)\leq s.
\de
So, $\tilde\phi_n(T_n)\in {\bf K}(s)$, which conflicts with (\ref{tildephi}). The proof is complete.
\end{proof}

Keeping the above lemma in mind, we prove upper bounds.

{\bf Proof of $(III)$.} First of all, for any $s\geq 0, \d>0, \t>0$, we take $\eta>s-\t>0$. So, by Lemma \ref{muees}, there exists a $r_1>0, \e_0>0$ such that for all $0<\e<\e_0$
\ce
\mu^\e(B^c(r_1)\cap\overline{\cD(A)})\leq \exp\(-\frac{\eta}{\e}\)\leq \exp\(-\frac{s-\t}{\e}\).
\de

Next, for $r_1$ and $R'$ in Lemma \ref{kbrks} corresponding to $\d, s$, by Lemma \ref{lam0t}, there exists a $\tilde T>0$ such that 
\be
\tilde\d:=\inf\{\Lambda_0^{\tilde T}(\phi): \phi\in D([0,\tilde T], \overline{\cD(A)}), |\phi(0)|\leq r_1, |\phi(\tilde T)|\geq R'\}>0.
\label{tild}
\ee
So, for any $n\in\mN$ set
\ce
\Gamma(r_1, s, \d, n):=\{\phi\in D([0,n\tilde T],\overline{\cD(A)}): |\phi(0)|\leq r_1, |\phi(i\tilde T)|\geq R', |\phi(i\tilde T)|\leq r_1, i=1, \cdots, n\},
\de
and there exists a $n'\in\mN$ such that 
\be
\inf\{\Lambda_0^{n'\tilde T}(\phi): \phi\in \Gamma(r_1, s, \d, n')\}>s.
\label{lambnp}
\ee
Indeed, if $\Lambda_0^{n'\tilde T}(\phi)=\infty$ for all $\phi\in \Gamma(r_1, s, \d, n')$, the left side of (\ref{lambnp}) is $\infty$ and (\ref{lambnp}) is right. If there exists a $\phi\in \Gamma(r_1, s, \d, n')$ such that $\Lambda_0^{n'\tilde T}(\phi)<\infty$, 
there exists a $(h,g)$ such that $\phi=X^{h,g}$ and 
\ce
\frac{1}{2}\int_0^{n'\tilde T}|h(t)|^2\dif t+\int_0^{n'\tilde T}\int_{\mU}\ell\(g(t,u)\)\nu(\dif u)\dif t<\infty.
\de
Note that
\be
&&\frac{1}{2}\int_0^{n'\tilde T}|h(t)|^2\dif t+\int_0^{n'\tilde T}\int_{\mU}\ell\(g(t,u)\)\nu(\dif u)\dif t\no\\
&=&\sum\limits_{i=1}^{n'}\left[\frac{1}{2}\int_{(i-1)\tilde T}^{i\tilde T}|h(t)|^2\dif t+\int_{(i-1)\tilde T}^{i\tilde T}\int_{\mU}\ell\(g(t,u)\)\nu(\dif u)\dif t\right]\no\\
&=&\sum\limits_{i=1}^{n'}\left[\frac{1}{2}\int_{0}^{\tilde T}|h_i(t)|^2\dif t+\int_{0}^{\tilde T}\int_{\mU}\ell\(g_i(t,u)\)\nu(\dif u)\dif t\right],
\label{npsum}
\ee
where $h_i(t)=h(t+(i-1)\tilde T), g_i(t,u)=g(t+(i-1)\tilde T,u)$ for $t\in[0,\tilde T]$ and $i=1,2,\cdots, n'$. Thus, $\phi_i(t):=\phi(t+(i-1)\tilde T)$ satisfies Eq.(\ref{xx0eq}) on $[0,\tilde T]$ with the initial value $\phi_i(0)=\phi((i-1)\tilde T)$ when $(h,g)$ is replaced by $(h_i, g_i)$. 
Since $|\phi_i(0)|\leq r_1, |\phi_i(\tilde T)|\geq R'$, by (\ref{tild}) it holds that $\Lambda_0^{\tilde T}(\phi_i)\geq \tilde\d$ and  
\ce
\frac{1}{2}\int_{0}^{\tilde T}|h_i(t)|^2\dif t+\int_{0}^{\tilde T}\int_{\mU}\ell\(g_i(t,u)\)\nu(\dif u)\dif t\geq \tilde\d,
\de
which together with (\ref{npsum}) implies that
\ce
\frac{1}{2}\int_0^{n'\tilde T}|h(t)|^2\dif t+\int_0^{n'\tilde T}\int_{\mU}\ell\(g(t,u)\)\nu(\dif u)\dif t\geq n'\tilde\d.
\de
We emphasize that the above inequality holds for all $(h,g)$ such that $\phi=X^{h,g}$. Therefore
\ce
\Lambda_0^{n'\tilde T}(\phi)\geq n'\tilde\d, \quad \forall \phi\in \Gamma(r_1, s, \d, n').
\de
Taking $n'\geq \frac{s}{\tilde\d}$, we have (\ref{lambnp}). 

Since $\Gamma(r_1, s, \d, n')$ is close in $D([0,n\tilde T],\overline{\cD(A)})$ and $B(r_1)\cap\overline{\cD(A)}$ is compact in $\overline{\cD(A)}$, (\ref{lambnp}) and Theorem \ref{dzuldpth} yield that for any $0<\e<\e_0$
\ce
\sup\limits_{x_0\in B(r_1)\cap\overline{\cD(A)}}\mP(X^\e(x_0)\in \Gamma(r_1, s, \d, n'))\leq \exp\left\{-\frac{s-\t}{\e}\right\}.
\de

Now, by combining the above deduction, the invariance of $\mu^\e$ implies that for any $T>0$
\be
&&\mu^\e(y\in\overline{\cD(A)}: \rho_{\overline{\cD(A)}}(y, {\bf K}(s))\geq \d)\no\\
&=&\int_{\overline{\cD(A)}}\mP(\rho_{\overline{\cD(A)}}(X^\e_{T}(x_0), {\bf K}(s))\geq\d)\mu^\e(\dif x_0)\no\\
&=&\int_{B^c(r_1)\cap\overline{\cD(A)}}\mP(\rho_{\overline{\cD(A)}}(X^\e_{T}(x_0), {\bf K}(s))\geq\d)\mu^\e(\dif x_0)\no\\
&&+\int_{B(r_1)\cap\overline{\cD(A)}}\mP(\rho_{\overline{\cD(A)}}(X^\e_{T}(x_0), {\bf K}(s))\geq\d)\mu^\e(\dif x_0)\no\\
&\leq&\int_{B^c(r_1)\cap\overline{\cD(A)}}\mP(\rho_{\overline{\cD(A)}}(X^\e_{T}(x_0), {\bf K}(s))\geq\d)\mu^\e(\dif x_0)\no\\
&&+\int_{B(r_1)\cap\overline{\cD(A)}} \mP(\rho_{\overline{\cD(A)}}(X^\e_{T}(x_0), {\bf K}(s))\geq\d, X^\e(x_0)\in \Gamma(r_1, s, \d, n'))\mu^\e(\dif x_0)\no\\
&&+\int_{B(r_1)\cap\overline{\cD(A)}}\mP(\rho_{\overline{\cD(A)}}(X^\e_{T}(x_0), {\bf K}(s))\geq\d, X^\e(x_0)\notin \Gamma(r_1, s, \d, n'))\mu^\e(\dif x_0)\no\\
&\leq&\mu^\e(B^c(r_1)\cap\overline{\cD(A)})+\sup\limits_{x_0\in B(r_1)\cap\overline{\cD(A)}}\mP(X^\e(x_0)\in \Gamma(r_1, s, \d, n'))\no\\
&&+\int_{B(r_1)\cap\overline{\cD(A)}}\mP(\rho_{\overline{\cD(A)}}(X^\e_{T}(x_0), {\bf K}(s))\geq\d, X^\e(x_0)\notin \Gamma(r_1, s, \d, n'))\mu^\e(\dif x_0)\no\\
&\leq& 2\exp\left\{-\frac{s-\t}{\e}\right\}\no\\
&&+\int_{B(r_1)\cap\overline{\cD(A)}}\mP(\rho_{\overline{\cD(A)}}(X^\e_{T}(x_0), {\bf K}(s))\geq\d, X^\e(x_0)\notin \Gamma(r_1, s, \d, n'))\mu^\e(\dif x_0).
\label{2st}
\ee

In the following, we estimate the last integration of the right side for the above inequality. Note that 
\be
&&\int_{B(r_1)\cap\overline{\cD(A)}}\mP(\rho_{\overline{\cD(A)}}(X^\e_{T}(x_0), {\bf K}(s))\geq\d, X^\e(x_0)\notin \Gamma(r_1, s, \d, n'))\mu^\e(\dif x_0)\no\\
&\leq&\int_{B(r_1)\cap\overline{\cD(A)}}\mP(\cup_{i=1}^{n'}\{\rho_{\overline{\cD(A)}}(X^\e_{T}(x_0), {\bf K}(s))\geq\d, X_{i\tilde T}^\e(x_0)\in B^c(r_1)\cap\overline{\cD(A)}\})\mu^\e(\dif x_0)\no\\
&&+\int_{B(r_1)\cap\overline{\cD(A)}}\mP(\cup_{i=1}^{n'}\{\rho_{\overline{\cD(A)}}(X^\e_{T}(x_0), {\bf K}(s))\geq\d, X_{i\tilde T}^\e(x_0)\in B(R')\cap\overline{\cD(A)}\})\mu^\e(\dif x_0)\no\\
&\leq&\int_{B(r_1)\cap\overline{\cD(A)}}\sum\limits_{i=1}^{n'}\mP\{\rho_{\overline{\cD(A)}}(X^\e_{T}(x_0), {\bf K}(s))\geq\d,  X_{i\tilde T}^\e(x_0)\in B^c(r_1)\cap\overline{\cD(A)}\}\mu^\e(\dif x_0)\no\\
&&+\int_{B(r_1)\cap\overline{\cD(A)}}\sum\limits_{i=1}^{n'}\mP\{\rho_{\overline{\cD(A)}}(X^\e_{T}(x_0), {\bf K}(s))\geq\d, X_{i\tilde T}^\e(x_0)\in B(R')\cap\overline{\cD(A)}\}\mu^\e(\dif x_0)\no\\
&=:&J_1+J_2.
\label{j1j2}
\ee
For $J_1$, from the invariance of $\mu^\e$ it follows that
\be
J_1&\leq&\sum\limits_{i=1}^{n'} \int_{B(r_1)\cap\overline{\cD(A)}}\mP\{X_{i\tilde T}^\e(x_0)\in B^c(r_1)\cap\overline{\cD(A)}\}\mu^\e(\dif x_0)\no\\
&\leq&\sum\limits_{i=1}^{n'} \int_{\overline{\cD(A)}}\mP\{X_{i\tilde T}^\e(x_0)\in B^c(r_1)\cap\overline{\cD(A)}\}\mu^\e(\dif x_0)\no\\
&=&\sum\limits_{i=1}^{n'} \mu^\e(B^c(r_1)\cap\overline{\cD(A)})\leq n'\exp\(-\frac{s-\t}{\e}\).
\label{j1}
\ee
For $J_2$, the Markov property of $X^\e(x_0)$ implies that for any $T>n'\tilde T$
\ce
J_2\leq\sum\limits_{i=1}^{n'}\sup\limits_{x_0\in B(R')\cap\overline{\cD(A)}}\mP\(\rho_{\overline{\cD(A)}}(X^\e_{T-i\tilde T}(x_0), {\bf K}(s))\geq\d\).
\de

Besides, by Lemma \ref{kbrks}, there exists a $T'>0$ such that for any $x_0\in B(R')\cap\overline{\cD(A)}$ and $T\geq T'$
\ce
\mP(\rho_{\overline{\cD(A)}}(X^\e_{T}(x_0), {\bf K}(s))\geq\d)&\leq& \mP(\rho_{D([0,T],\overline{\cD(A)})}(X^\e(x_0), {\bf K}_{B(R'), 0, T}(s))\geq\d/2)\\
&\leq& \mP(\rho_{D([0,T],\overline{\cD(A)})}(X^\e(x_0), {\bf K}_{\{x_0\}, 0, T}(s))\geq\d/2).
\de
Note that ${\bf K}_{\{x_0\}, 0, T}(s)=\Phi_{x_0}(s)$. Thus, Remark \ref{fwthagst} $(iii)$ yields that for any $x_0\in B(R')\cap\overline{\cD(A)}$
\ce
\mP(\rho_{D([0,T],\overline{\cD(A)})}(X^\e(x_0), {\bf K}_{\{x_0\}, 0, T}(s))\geq\d/2)\leq \exp\left\{-\frac{s-\t}{\e}\right\},
\de
and
\ce
\mP(\rho_{\overline{\cD(A)}}(X^\e_{T}(x_0), {\bf K}(s))\geq\d)\leq \exp\left\{-\frac{s-\t}{\e}\right\}.
\de

Finally, we take $T=T'+n'\tilde T$ and obtain that
\ce
J_2\leq n'\exp\left\{-\frac{s-\t}{\e}\right\},
\de
which together with (\ref{2st})-(\ref{j1}) implies that
\ce
\mu^\e(y\in\overline{\cD(A)}: \rho_{\overline{\cD(A)}}(y, {\bf K}(s))\geq \d)\leq (2+2n')\exp\left\{-\frac{s-\t}{\e}\right\}.
\de
Thus, the proof is complete.


\begin{thebibliography}{999}

\bibitem{by} J. Bao and C. Yuan: Large deviations for neutral functional SDEs with jumps, {\it Stochastics}, 87(2015)48-70.

\bibitem{bfz} R. Bai, C. Feng and H. Zhao: Large deviations principle for invariant measures of stochastic Burgers equations, http://arxiv.org/abs/2409.14234.

\bibitem{bc} Z. Brzezniak and S. Cerrai: Large deviations principle for the invariant measures of the 2D stochastic Navier-Stokes equations on a torus, {\it J. Funct. Anal.}, 273(2017)1891-1930.

\bibitem{bcd} A. Budhiraja, J. Chen and P. Dupuis: Large deviations for stochastic partial differential equations driven by Poisson random measure, {\it Stoch. Process. Appl.}, 123 (2013)523-560.

\bibitem{bdv1} A. Budhiraja, P. Dupuis and V. Maroulas: Large deviations for infinite dimensional stochastic dynamical systems, {\it The Annals of Probability}, 36(2008)1390-1420.

\bibitem{bdv2} A. Budhiraja, P. Dupuis and V. Maroulas: Variational representations for continuous time processes, {\it Ann. Inst. Henri Poincar\'e Probab. Statist.}, 47(2011)725-747.

\bibitem{cr1} S. Cerrai and M. R\"ockner: Large deviations for invariant measures of stochastic reaction-diffusion systems with multiplicative noise and non-Lipschitz reaction term, {\it Ann. Inst. H. Poincare Probab. Statist.}, 41(2005)69-105.

\bibitem{cr2} S. Cerrai and M. R\"ockner: Large deviations for stochastic reaction-diffusion systems with multiplicative noise and non-Lipschitz reaction term, {\it The Annals of Probability}, 32(2004)1100-1139.

\bibitem{cp} S. Cerrai and N. Paskal: Large deviations principle for the invariant measures of the 2D stochastic Navier-Stokes equations with vanishing noise correlation, {\it Stoch PDE: Anal. Comp.}, 10(2022)1651-1681.

\bibitem{cdjz} L. Chen, Z. Dong, J. Jiang and J. Zhai: On limiting behavior of stationary measures for stochastic evolution systems with small noise intensity, {\it Sci. China. Math.}, 63(2020)1463-1504.

\bibitem{cm} F. Chenal and A. Millet: Uniform large deviations for parabolic SPDEs and applications, {\it Stochastic Processes and their Applications}, 72(1997)161-186.

\bibitem{dz} A. Dembo and O. Zeitouni: {\it Large deviations techniques and applications}, Springer Science \& Business Media, 2009.

\bibitem{ff} C. A. D. Fraga Filho: On the boundary conditions in Lagrangian particle methods and the physical foundations of continuum mechanics, {\it Continuum Mech. Thermodyn.}, 31(2019)475-489.

\bibitem{fw} M. I. Freidlin and A. D. Wentzell: {\it Random Perturbations of Dynamical Systems}, Springer, New York, Third Edition, 2012.

\bibitem{eg} E. Gautier: Uniform large deviations for the nonlinear Schr\"odinger equation with multiplicative noise, {\it Stochastic Processes and their Applications}, 115(2005)1904-1927.

\bibitem{gw} Y. Guan and J. Wu: Exponential ergodicity for non-Lipschitz multivalued stochastic differential equations with L\'evy jumps, {\it Infinite Dimensional Analysis, Quantum Probability and Related Topics}, 20(2017)1750002(21pages).

\bibitem{hjh} E. M. Hanks, D. S. Johnson and M. B. Hooten: Reflected Stochastic Differential Equation Models for Constrained Animal Movement, {\it Journal of Agricultural, Biological, and Environmental Statistics}, 22(2017)353-372.

\bibitem{hwy} S. He, J. Wang and J. Yan: {\it Semimartingale Theory and Stochastic Calculus}, Science Press and CRC Press Inc., 1992.

\bibitem{hjly} W. Huang, M. Ji, Z. Liu and Y. Yi: Concentration and limit behaviors of stationary measures, {Physica D}, 369(2018)1-17.

\bibitem{jwzz} J. Jiang, J. Wang, J. Zhai and T. Zhang: Uniform large deviations and metastability of random dynamical systems, http://arxiv.org/abs/2402.16522.

\bibitem{km} A. Kumar and M. T. Mohan: Uniform large deviation principle for the solutions of two-dimensional stochastic Navier-Stokes equations in vorticity form, {\it Applied Mathematics \& Optimization},  90(2024)9.

\bibitem{lw} S. C. Leite and R. J. Williams: A constrained Langevin approximation for chemical reaction networks, {\it Ann. Appl. Probab.}, 2019, 29: 1541-1608.

\bibitem{xmfx} X. Ma and F. Xi: Large deviations for invariant measures of stochastic differential equations with jumps, {\it Stochastics}, 91(2019)528-552.

\bibitem{Mar} V. Maroulas: Uniform large deviations for infinite dimensional stochastic systems with jumps, {\it Mathematika}, 57(2011)175-192.

\bibitem{dm1} D. Martirosyan: Large deviations for invariant measures of the white-forced $2D$ Navier-Stokes equation, {\it Journal of Evolution Equations}, 18(2018)1245-1265.

\bibitem{dm2} D. Martirosyan: Large deviations for stationary measures of stochastic nonlinear wave equations with smooth white noise, {\it Comm. Pure Appl. Math.}, 70(2017)1754-1797.

\bibitem{mr} J. L. Menaldi and M. Robin: Reflected diffusion processes with jumps, {\it Ann. Probab.}, 13(1985)319-341.

\bibitem{q1} H. Qiao: Asymptotic behaviors of multiscale multivalued stochastic systems with small noises, http://arxiv.org/abs/2306.06922.

\bibitem{q2} H. Qiao: Exponential ergodicity for SDEs with jumps and non-Lipschitz coefficients, {\it J.Theor. Probab.}, 27(2014)137-152.

\bibitem{qg} H. Qiao and J. Gong: Backward multivalued McKean-Vlasov SDEs and associated variational inequalities, {\it Discrete and Continuous Dynamical Systems-S}, 16(2023)819-845. 

\bibitem{zq} Z. Qiu: The ergodicity and uniform large deviations for the 1D stochastic Landau-Lifshitz-Bloch equation, {\it Stochastic Analysis and Applications}, 42(2024)963-985.

\bibitem{rxz} J. Ren, S. Xu and X. Zhang: Large deviations for multivalued stochastic differential equations, {\it J. Theor. Probab.}, 23(2010)1142-1156.

\bibitem{rwz} J. Ren, J. Wu and H. Zhang: General large deviations and functional iterated logarithm law for multivalued stochastic differential equations, {\it J. Theor. Probab.}, 28(2015)550-586.

\bibitem{ms} M. Salins: Equivalences and counterexamples between several definitions of the uniform large deviations principle, {\it Probab. Surv.}, 16(2019)99-142.

\bibitem{ms1} M. Salins: Systems of small-noise stochastic reaction-diffusion equations satisfy a large deviations principle that is uniform over all initial data, {\it Stochastic Processes and their Applications}, 142(2021)159-194.

\bibitem{sbd} M. Salins, A. Budhiraja and P. Dupuis: Uniform large deviation principles for Banach space valued stochastic evolution equations, {\it Transactions of the AMS}, 372(2019)8363-8421.

\bibitem{ss} M. Salins and L. Setayeshgar: Uniform large deviations for a class of Burgers-type stochastic partial differential equations in any space dimension, {\it Potential Analysis}, 58(2023)181-201.

\bibitem{so} R. Sowers: Large deviations for the invariant measures of a reaction-diffusion equation with non Gaussian perturbations, {\it Probab. Theory Relat. Fields}, 92(1992)393-421.

\bibitem{wangiv1} B. Wang: Large deviations of invariant measures of stochastic reaction-diffusion equations on unbounded domains, {\it Journal of Statistical Physics}, 191(2024)96.

\bibitem{wangiv2} B. Wang: Large deviation principles of invariant measures of stochastic reaction-diffusion lattice systems, http://arxiv.org/abs/2405.02720.

\bibitem{wang1} B. Wang: Uniform large deviation principles of fractional reaction-diffusion equations driven by superlinear multiplicative noise on $\mR^n$, http://arxiv.org/abs/2406.08722.

\bibitem{wang2} B. Wang: Uniform large deviations of fractional stochastic equations with polynomial drift on unbounded domains, {\it Stochastics and Dynamics}, 23(2023)2350049.

\bibitem{wy} J. Wang and H. Yang: Uniform large deviation principles for SDEs under locally weak monotonicity conditions, http://arxiv.org/abs/2409.02153.

\bibitem{jw} J. Wu: Uniform large deviations for multivalued stochastic differential equations with Poisson jumps, {\it Kyoto J. Math.}, 51(2011)535-559.

\bibitem{hz} H. Zhang: Large deviations for invariant measures of multivalued stochastic differential equations, {\it Stochastic Analysis and Applications}, 40(2022)798-811.

\bibitem{tz} T. Zhang: Large deviations for invariant measures of SPDEs with two reflecting walls, {\it Stochastic Processes and their Applications}, 122(2012)3425-3444.

\bibitem{zl} Y. Zhang and X. Li: Large deviation principles of invariant measures of stochastic discrete wave equations with multiplicative noise and nonlinear damping, {\it Bull. Sci. math.}, 204(2025)103670.

\end{thebibliography}
\end{document}